\documentclass[10pt, a4]{amsart}
\usepackage{amsmath, amsfonts, amssymb,  mathrsfs, color, amsthm, epsfig, color, graphics, old-arrows}
  \usepackage[all]{xy}
 \usepackage[latin1]{inputenc}
\usepackage{hyperref}

\newtheorem{definition}{\bf Definition}[section]
\newtheorem{prop}[definition]{{\bf Proposition}}
\newtheorem{rem}[definition]{{\bf Remark}}
\newtheorem{lemma}[definition]{{\bf Lemma}}
\newtheorem{cor}[definition]{{\bf Corollary}}
\newtheorem{hypothesis}[definition]{{\bf Hypothesis}}
\newtheorem{theorem}[definition]{{\bf Theorem}}

\DeclareTextFontCommand{\textcyr}{\cyr}

\usepackage{anysize} 
\marginsize{2.5cm}{2.5cm}{2cm}{2cm}

\DeclareMathOperator{\tr}{T}

\DeclareMathOperator{\sumperp}{\stackrel{\perp}{\oplus}}

\DeclareMathOperator*{\argmin}{argmin} 

\makeatletter
\renewcommand{\theequation}{\ifnum \c@section>\z@ \thesection.\fi
	\@arabic\c@equation}
\@addtoreset{equation}{section}
\makeatother

\newcommand{\F}{\mathcal F}
\newcommand{\G}{\mathcal G}
\newcommand{\K}{{\mathcal K}}
\newcommand{\B}{\mathcal O}

\newcommand{\Feps}{{\mathcal F_\varepsilon}}
\newcommand{\Beps}{\mathcal O_\varepsilon}
\newcommand{\dx}{\,{\rm d}x}
\newcommand{\ds}{\,{\rm d}s}

\newcommand{\forallt}{\qquad\text{for all }}

\renewcommand{\O}{\mathcal O}

\title[Asymptotic limit for the Stokes and Navier--Stokes problems]{Asymptotic limit for the Stokes and Navier--Stokes problems in a planar  domain with a vanishing hole}

\author{Alexandre Munnier}
\date{\today}

\begin{document}
\maketitle
%
\begin{abstract}
We show that the eigenvalues  of the Stokes operator in a domain with a small hole
converge to the eigenvalues 
of the Stokes operator in the whole domain, when the diameter of the hole tends to 0. The convergence of the eigenspaces 
and the  convergence of the Stokes semigroup are also established. Concerning the Navier--Stokes equations, we prove that  
the vorticity of the solution in the perforated domain  converges as the hole shrinks to a point $r$ to the vorticity of the solution in the punctured
domain (i.e. the whole domain with the point $r$ removed).
The main ingredients of the 
analysis are a suitable decomposition of the vorticity space, the formalism elaborated in \cite{Lequeurre:2020aa} and some basics of potential 
theory. 

\end{abstract}
 %
 \section{Introduction}
 Let $\F$ be an open, bounded and smooth domain in $\mathbb R^2$.
The  Lebesgue space $\mathbf L^2(\F)=L^2(\F,\mathbb R^2)$ 
 and the Sobolev space $\mathbf H^1_0(\F)=H^1_0(\F,\mathbb R^2)$ are equipped with their usual scalar products and 
 the Hilbert spaces:
\begin{equation}
\label{def:J0}
\mathbf J_0(\F)=\big\{u\in \mathbf L^2(\F)\,:\, \nabla\cdot u=0\text{ in }\F,\,u\cdot n=0\text{ on }\partial\F\big\}
\qquad\text{and}\qquad
\mathbf J_1(\F)=\mathbf J_0(\F)\cap \mathbf H^1_0(\F),
\end{equation}
 are provided respectively with the scalar products:
\begin{subequations}
\begin{alignat}{3}
\label{def:scal_J0}
(u,v)_{\mathbf J_0(\F)}&= (u ,v)_{\mathbf L^2(\F)} &&\forallt u,v\in\mathbf J_0(\F),\\
\label{def:scal_J1}
(u,v)_{\mathbf J_1(\F)}&=(\nabla^\perp\cdot u,\nabla^\perp\cdot v)_{L^2(\F)}&&\forallt u,v\in\mathbf J_1(\F).
\end{alignat}
\end{subequations}
 In \eqref{def:J0},  $n$ stands for the unit outer normal vector to $\partial\F$. In \eqref{def:scal_J1}   and subsequently in the paper, for every $x=(x_1,x_2)\in\mathbb R^2$, the notation $x^\perp$ is used to represent the vector $(-x_2,x_1)$. Identifying $\mathbf J_0(\F)$ with its dual space $\mathbf J_0'(\F)$ and denoting by $\mathbf J_{-1}(\F)$ the dual space of $\mathbf J_1(\F)$ with respect to the pivot  $\mathbf J_0(\F)$, we obtain a Gelfand triple of Hilbert spaces: 
$$\mathbf J_1(\F)\subset \mathbf J_0(\F)\subset \mathbf J_{-1}(\F),$$
both inclusions being continuous and dense. The Stokes operator $\mathsf J_\F$ is the unbounded operator on $\mathbf J_0(\F)$ of domain
$$D(\mathsf J_\F)=\big\{u\in\mathbf J_1(\F)\,:\, (u,\cdot)_{\mathbf J_1(\F)}\in\mathbf J_0'(\F)\big\},$$
and defined for every $u\in D(\mathsf J_\F)$ by means of the Riesz representation Theorem by:
$$(\mathsf J_\F u,\cdot)_{\mathbf J_0(\F)}=(u,\cdot)_{\mathbf J_1(\F)}.$$
The spectrum 
of $\mathsf J_\F$ consists in a sequence of monotonically ordered positive eigenvalues  $(\lambda_k^\F)_{k\geqslant 1}$ that tends to $+\infty$. The eigenvalues are counted with their multiplicity. They meet the Courant--Fischer min-max principle:
\begin{equation}
\label{eq:rayleigh}
\lambda_k^\F=\min_{V\in W_k^J}\max_{{\theta\in V}\atop{\theta\neq 0}}\frac{\|\theta\|_{\mathbf J_1(\F)}^2}{\|\theta\|_{\mathbf J_0(\F)}^2}\forallt k\in\mathbb N,\,k\geqslant 1,
\end{equation}
where $W_k^J$ stands for the set of all the subspaces of dimension $k$ in $\mathbf J_1(\F)$. 
For every positive integer $k$, we denote by $I_k^\F$ the set of all the indices $j$ such that $\lambda_j^\F=\lambda_k^\F$ and
 by $\varLambda_k^\F$ the eigenspace associated with the eigenvalue $\lambda_k^\F$. This implies in particular that $\dim \varLambda_k^\F=\# I_k^\F$ 
 and that $\varLambda_j^\F=\varLambda_k^\F$ if $j\in I_k^\F$. To every eigenvalue (still counted with algebraic multiplicity) we associate an eigenfunction $u_k^\F$ chosen in such a way that 
 the family $\big\{u_k^\F,\,k\geqslant 1\big\}$ is a Riesz Basis orthonormal in $\mathbf J_0(\F)$ and orthogonal in $\mathbf J_1(\F)$.
%
%
\par
Our first purpose is to  study the behavior of the spectrum of $\mathsf J_\F$ when the domain $\F$ has a small hole whose diameter tends to 0.
More precisely, the geometry considered is as follows: $\G$ is an open, bounded and simply connected domain with a smooth boundary denoted by 
$\Gamma$. For every $\varepsilon$ (a real number ranging from 0 to some positive real number $\varepsilon_0$) we define a smooth 
and simply connected domain 
$\Beps$ (subsequently referred to as the ``hole'' in the domain or the ``obstacle'' in the fluid). The boundary of $\Beps$
is denoted by $\Sigma_\varepsilon$  and we assume that  there exists a point $r\in\G$ such that:
$$\overline{\mathcal D(r,\varepsilon_0)}\subset \G
\quad\text{and}\quad\overline{\Beps}\subset \mathcal D(r,\varepsilon)\forallt \varepsilon\in(0,\varepsilon_0),$$
where $\mathcal D(r,\varepsilon)$ stands for the disk of center $r$ and radius $\varepsilon$. As $\varepsilon$ tends to 0, we shall write that $\Beps$ 
``shrinks'' (or vanishes) to a point. The perforated domain (occupied by the fluid) is 
$\Feps=\G\setminus\overline{\B}_\varepsilon$  and thereby its boundary is the 
disjoint union $\Gamma\cup\Sigma_\varepsilon$. Notice that the notion of ``vanishing'' obstacle as defined here is more general 
than the one considered in \cite{Iftimie:2006aa} or \cite{Lacave:2017aa} for instance, where the obstacle is the homothetic image of a reference configuration. 
Our definition is roughly the same as in \cite{He:2019aa} (where the more intricated case of moving obstacles is addressed).
\par
The statement of the first main result yet requires making precise the notion of subspace convergence: 
For every $\varepsilon\in(0,\varepsilon_0)$, let $W_\varepsilon$ be a closed subspace in a Hilbert space $H$ and let 
the orthogonal projection on  $W_\varepsilon$
be denoted by $\Pi_{W_\varepsilon}$. Let $W$ be another closed subspace of $H$ and let 
the orthogonal projection on $W$ be denoted by $\Pi_W$. We shall write  that $W_\varepsilon\longrightarrow W$ as $\varepsilon\longrightarrow 0$ when:
$$\sup_{\theta\in H\atop \theta\neq 0}\frac{\|\Pi_W\theta -\Pi_{W_\varepsilon} \theta\|_H}{\|\theta\|_H}
\longrightarrow 0\quad\text{as }\varepsilon\longrightarrow 0.$$
Extending the functions by $0$ inside $\Beps$, we can assume that for every $\varepsilon\in(0,\varepsilon_0)$, 
$\mathbf J_0(\Feps)$ is a closed subspace of $\mathbf J_0(\G)$ and that $\mathbf J_1(\Feps)$ is a closed subspace of 
$\mathbf J_1(\G)$. In the same manner, the eigenspaces $\varLambda_k^{\Feps}$ can be considered as closed subspaces of $\mathbf J_1(\G)$.
\begin{theorem}[Convergence of eigenvalues and eigenspaces]
\label{main:theo}
Assume that $\Beps$ shrinks to a point as $\varepsilon$ goes to $0$. Then, for every positive integer $k$:
\begin{subequations}
\label{conv_eigen:gene}
\begin{alignat}{3}
\label{conv_eigen}
\lambda_k^{\Feps}&\longrightarrow \lambda_k^\G&\quad&\text{as }\varepsilon\longrightarrow 0,\\
\label{conv_eigen:2}
\bigoplus_{j\in I_k^\G}\varLambda_j^{\Feps}&\longrightarrow\varLambda_k^\G\quad\text{in}\quad\mathbf J_0(\G)&\quad&\text{as }\varepsilon\longrightarrow 0.
\end{alignat}
\end{subequations}
\end{theorem}
%
We emphasize that in \eqref{conv_eigen:2} the sum ranges over all the indices $j$ such that $\lambda_j^\G=\lambda_k^\G$ (because 
some eigenvalues can be different when $\varepsilon>0$ and may eventually meet when $\varepsilon=0$).
\par
%
The asymptotic behavior of  solutions of the (stationary) Stokes equations in a domain with a small hole has been widely investigated; see for instance 
\cite{Caubet:2016aa} and references therein. However, the asymptotic limit of the eigenvalues and eigenspaces of the Stokes operator has 
not been dealt with so far. 
\par
The proof of Theorem~\ref{main:theo} rests on a restatement of the Stokes operator in term of   so--called non--primitive variables (stream function and vorticity). This task 
was carried out in the paper 
\cite{Lequeurre:2020aa} and briefly summarized  later on. Then, the conclusion of the theorem derives from a suitable
decomposition of the vorticity space (established 
in Section~\ref{SEC:decomp_vorti}) for a perforated domain.
\par
%
The convergence results \eqref{conv_eigen:gene} allow quite easily deriving convergence results for the Stokes semigroup.
Denote by $\big\{\mathsf T_\G(t),\,t\geqslant 0\big\}$ the Stokes semigroup whose infinitesimal generator is $\mathsf J_\G$ (the Stokes operator 
for the domain $\G$). For every 
$\theta\in \mathbf J_0(\G)$, we have the classical expression:
\begin{equation}
\label{def_semigroup}
\mathsf T_\G(t)\theta=\sum_{j\geqslant 1}\big(\theta,u_j^\G\big)_{\mathbf J_0(\G)}e^{-\lambda_j^\G t}u_j^\G,\qquad t\geqslant 0.
\end{equation}
In the same manner and for every $\varepsilon\in(0,\varepsilon_0)$, we can define for the domain $\Feps$ the semigroup $\big\{\mathsf T_{\Feps}(t),\,t\geqslant 0\big\}$ whose 
infinitesimal generator is the Stokes operator $\mathsf J_{\Feps}$. Thus:
\begin{equation}
\label{def_semigroup}
\mathsf T_\Feps(t)\theta=\sum_{j\geqslant 1}\big(\theta,u_j^\Feps\big)_{\mathbf J_0(\G)}e^{-\lambda_j^\Feps t}u_j^\Feps,\qquad t\geqslant 0,
\end{equation}
and this expression makes sense for every $\theta\in \mathbf J_0(\G)$.
\begin{cor}
\label{cor:1}
The following limit holds for every $T>0$ and every $\theta\in \mathbf J_0(\G)$:
\begin{subequations}
\label{conv_eigen-spaces}
\begin{equation}
\label{conv_eigen-spaces_1}
\sup_{t\in[0,T]}e^{\lambda_1^\G t}{\big\| \mathsf T_\G(t)\theta-\mathsf T_{\Feps}(t)\theta \big\|_{\mathbf J_0(\G)}}
\longrightarrow 0\quad\text{as}
\quad\varepsilon\longrightarrow 0.
\end{equation}
For every compact set $K\subset \mathbf J_0(\G)$ and every $T>0$:
\begin{equation}
\label{conv_eigen-spaces_2}
\sup_{t\in[0,T]\atop 
\theta\in K,\,\theta\neq 0}
e^{\lambda_1^\G t}\frac{\big\| \mathsf T_\G(t)\theta-\mathsf T_{\Feps}(t)\theta \big\|_{\mathbf J_0(\G)}}{\|\theta\|_{\mathbf J_0(\G)}}\longrightarrow 0\quad\text{as}
\quad\varepsilon\longrightarrow 0.
\end{equation}
\end{subequations}
\end{cor}
The exponential decay property:
$$\|\mathsf T_\G(t)\theta\|_{\mathbf J_0(\G)}\leqslant \|\theta\|_{\mathbf J_0(\G)}e^{-\lambda_1^\G t}
\forallt \theta\in\mathbf J_0(\G),$$ 
explains the role played by the term $e^{\lambda_1^\G t}$ in the estimate \eqref{conv_eigen-spaces}. 
\par
We turn now our attention to the Navier--Stokes equations. 
For every $\varepsilon\in(0,\varepsilon_0)$, let $u_\varepsilon^0$ be given in $\mathbf J_0(\Feps)$ and let $u_\varepsilon$ be the unique 
 function in
$$L^2(\mathbb R_+;\mathbf J_1(\Feps))\cap C(\mathbb R_+;\mathbf J_0(\Feps))
\cap H^1(\mathbb R_+;\mathbf J_{-1}(\Feps)),$$
that solves the following Cauchy problem for every $\theta\in \mathbf J_1(\Feps)$:
\begin{subequations}
\label{NS_classical}
\begin{alignat}{3}
\label{main_cauchy2}
\frac{\rm d}{{\rm d}t}(u_\varepsilon,\theta)_{\mathbf J_0(\Feps)}+\nu (u_\varepsilon,\theta)_{\mathbf J_1(\Feps)}-
((u_\varepsilon\cdot\nabla)\theta, u_\varepsilon)_{\mathbf L^2(\Feps)}&=0&\quad&\text{on }\mathbb R_+\\
u_\varepsilon(0)&=u^0_\varepsilon&&\text{in }\Feps.
\end{alignat}
\end{subequations}
%
\begin{theorem}
\label{main_theo_2}
 Assume that there exists $u^0\in \mathbf J_0(\G)$ such that 
$u_\varepsilon^0\longrightharpoonup u^0$ weak in $\mathbf J_0(\G)$ (here and subsequently, $u_\varepsilon^0$ and $u_\varepsilon$ 
are extended by $0$ inside $\Beps$). Then, as $\varepsilon$ goes to $0$:
\begin{subequations}
\label{convergence_ueps}
\begin{alignat}{3}
u_\varepsilon&\longrightharpoonup u&\text{ weak--}\star&\text{ in }L^\infty(\mathbb R_+;\mathbf J_0(\G)),\\
u_\varepsilon&\longrightarrow u&\text{ strong}&\text{ in }L^2_{\ell{\rm oc}}(\mathbb R_+;\mathbf J_0(\G)),\\
u_\varepsilon&\longrightharpoonup u&\text{ weak}&\text{ in }L^2(\mathbb R_+;\mathbf J_1(\G)),
\end{alignat}
\end{subequations}
where the function $u$ belongs to:
$$L^2(\mathbb R_+;\mathbf J_1(\G))\cap C(\mathbb R_+;\mathbf J_0(\G))
\cap H^1(\mathbb R_+;\mathbf J_{-1}(\G)),$$
and solves the Cauchy problem for every $\theta\in \mathbf J_1(\G)$:
\begin{subequations}
\label{convergeence_omegaeps}
\begin{alignat}{3}
\label{main_cauchy3}
\frac{\rm d}{{\rm d}t}(u,\theta)_{\mathbf J_0(\G)}+\nu (u,\theta)_{\mathbf J_1(\G)}-
((u\cdot\nabla)\theta, u)_{\mathbf L^2(\G)}&=0&\quad&\text{on }\mathbb R_+\\
u(0)&=u^0&&\text{in }\G.
\end{alignat}
\end{subequations}
Let now $\chi$ be in $\mathscr D(\G\setminus\{r\})$ and denote by $\omega^0_\varepsilon$ the vorticity of the velocity field $u_\varepsilon^0$.  Assume that for every $\varepsilon\in(0,\varepsilon_0)$, $\omega_\varepsilon^0$ is in $V_0(\Feps)$ 
and that the quantity $\|\chi\omega^0_\varepsilon\|_{L^2(\G)}$ is uniformly bounded. Then, 
as $\varepsilon$ goes to $0$:
\begin{subequations}
\label{gbhnj}
\begin{alignat}{3}
\chi\omega_\varepsilon&\longrightharpoonup \chi\omega&\text{ weak--}\star&\text{ in }L^\infty(\mathbb R_+;L^2(\G)),\\
\label{conv_palenstro}
\nabla(\chi\omega_\varepsilon)&\longrightharpoonup \nabla(\chi\omega)&\text{ weak }&\text{ in }
L^2(\mathbb R_+;\mathbf L^2(\G)),
\end{alignat}
\end{subequations}
where $\omega=\nabla^\perp\cdot u$ and $\omega_\varepsilon=\nabla^\perp\cdot u_\varepsilon$. 
\end{theorem}

Although stated in a different and more intricate  contexte (exterior domain or moving obstacles), the convergence results \eqref{convergence_ueps} 
meet those obtained in \cite{Iftimie:2006aa} and \cite{He:2019aa} and to this extent cannot be considered as new. However, we shall provide a completely different and more simple proof based 
on the stream--vorticity formulation of the Navier--Stokes equations introduced in \cite{Lequeurre:2020aa} and involving a different 
compactness argument. In contrast, the convergence \eqref{gbhnj} of the vorticity is new.  

As already mentioned  in \cite{Lequeurre:2020aa}, 
the analysis of the solutions to the Navier--Stokes equations in a planar domain is tightly related to the analysis of the harmonic functions and more precisely 
on some $L^2$--mass concentration properties near the boundaries of the domain. This provides an efficient and original strategy 
to deal with the problem.
\par
The rest of the paper is organized as follows: The next Section is a short summary of results from \cite{Lequeurre:2020aa}. 
Section~\ref{SEC:decomp_vorti} is dedicated to technical lemmas addressing mainly $L^2$--mass concentration properties of harmonic functions 
in a perforated domain. A theorem describing the structure of the vorticity space is also provided. The proofs of Theorem~ \ref{main:theo}
 and Corollary~\ref{cor:1} are carried out in Section~\ref{SEC:proof_of_theo1} and the last section contains the proof of 
 Theorem~\ref{main_theo_2}.
\section{The  Navier--Stokes equations in non--primitive variables}
Let us put aside for a while the perforated domain   $\Feps$ and  consider back as in the Introduction
the more general domain  simply denoted by $\F$.  In addition of being smooth and bounded, the domain  
$\F$ is  also assumed to be $N$--connected  ($N$ a nonnegative integer). The boundary  of $\F$
can be split into a disjoint union of smooth Jordan curves:
\begin{equation}
\partial\F=\Big(\bigcup_{k=1}^N \Sigma_k\Big)\cup\Gamma.
\end{equation}
The curves $\Sigma_k$ for $k\in\{1,\ldots,N\}$ are the inner boundaries of $\F$ while $\Gamma$ is the outer boundary. The Hilbert spaces:
\begin{subequations}
\label{def_S_spaces}
\begin{align}
S_0(\F)&=\{\psi\in  H^1(\F)\,:\,\, \psi|_{\Gamma}=0\quad\text{and}\quad \psi|_{\Sigma_j}=c_j,\quad c_j\in\mathbb R,\quad j=1,\ldots,N\},\\
S_1(\F)&=\bigg\{\psi\in S_0(\F)\cap H^2(\F)\,:\,\frac{\partial \psi}{\partial n}\Big|_{\partial\F}=0\bigg\},
\end{align}
%
are provided with the scalar products:
\begin{align}
(\psi_1,\psi_2)_{S_0(\F)}&=(\nabla \psi_1,\nabla\psi_2)_{\mathbf L^2(\F)}\forallt \psi_1,\psi_2\in  S_0(\F),\\
(\psi_1,\psi_2)_{S_1(\F)}&=(\Delta\psi_1,\Delta\psi_2)_{L^2(\F)}\forallt \psi_1,\psi_2\in  S_1(\F).
\end{align}
\end{subequations}
The space $S_1(\F)$ is continuously and densely embedded in $S_0(\F)$. Using the latter as pivot space and denoting by $S_{-1}(\F)$ 
the dual  of $S_1(\F)$, we obtain a Gelfand triple of Hilbert spaces:
$$S_1(\F)\subset S_0(\F)\subset S_{-1}(\F).$$
This provides the suitable functional framework to deal with the Navier--Stokes equations in stream function formulation. Thus:
%
\begin{theorem}[Well posedness of the weak NS equations in stream function formulation]
\label{theo_weak_NS}
For any $\psi^0\in S_0(\F)$, there exists a unique function:
$$\psi\in H^1(\mathbb R_+;S_{-1}(\F))\cap C(\mathbb R_+;S_0(\F))\cap L^2(\mathbb R_+;S_1(\F)),$$
satisfying for every $\theta\in S_1(\F)$ the Cauchy problem:
\begin{subequations}
\label{varia4}
\begin{alignat}{3}
\frac{\rm d}{{\rm d}t}(\psi,\theta)_{S_0(\F)}+\nu(\psi,\theta)_{S_1(\F)}-(D^2\theta \nabla^\perp\psi,\nabla\psi)_{\mathbf L^2(\F)}
&=0&\quad&\text{on }\mathbb R_+,\\
\psi(0)&=\psi^0&&\text{in }\F,
\end{alignat}
\end{subequations}
where $D^2\theta$ is the Hessian tensor field of $\theta$ in $\G$. 
\end{theorem}
The proof of this result (as the proofs of all the results stated in this section) can be found 
in \cite{Lequeurre:2020aa}.
\par
%
We shall now established the expression of the Navier--Stokes equations in vorticity formulation. Let
  $\mathfrak H(\F)$ stand for the closed space of the harmonic functions $h$ in $L^2(\F)$ verifying, for every smooth Jordan curve $\mathscr C$ 
included in $\F$:
\begin{equation}
\label{nul_flux}
\int_{\mathscr C}\frac{\partial h}{\partial n}\ds=0,
\end{equation}
where $n$ is the unit normal vector to the curve $\mathscr C$ (when $\F$ is simply connected, this condition is automatically  
satisfied by any harmonic function in 
$L^2(\F)$).
 In the flux condition \eqref{nul_flux}, the normal derivative of $h$ on $\mathscr C$ is well defined as an element of $H^{-1/2}(\mathscr C)$
 (so the integral should be understood as a duality bracket). Next, we introduce $V_0(\F)=\mathfrak H(\F)^\perp$ so that:
\begin{equation}
\label{decomp_L2}
L^2(\F)=V_0(\F)\sumperp \mathfrak H(\F).
\end{equation}
 %

 We recall that (see \cite{Guermond:1994aa}, \cite{Lequeurre:2020aa} for a  proof):
 \begin{prop}
 \label{prop:iso}
The following operators are isometries:
\begin{subequations}
\label{isom:1}
\begin{alignat}{3}
\label{isom:1-1}
\nabla^\perp:S_0(\F)&\longrightarrow\mathbf J_0(\F),&\qquad&\psi\longmapsto \nabla^\perp\psi,\\
\label{isom:1-2}
\nabla^\perp:S_1(\F)&\longrightarrow\mathbf J_1(\F),&\qquad&\psi\longmapsto \nabla^\perp\psi,\\
\label{isom:1-3}
\Delta:S_1(\F)&\longrightarrow V_0(\F),&\qquad&\psi\longmapsto \Delta\psi.
\end{alignat}
\end{subequations}
The inverse of the last operator is usually referred to as the Biot-Savart operator.
\end{prop}
%
For a velocity field $u$ in the space $\mathbf J_1(\F)$, the function $\psi$ such that $\nabla^\perp\psi=u$ is the associated stream function and $\omega=\Delta\psi=\nabla^\perp\cdot u$ 
is the vorticity field. From Proposition~\ref{prop:iso}, it can easily be deduced that the formulation \eqref{varia4} of the Navier--Stokes 
equations is equivalent to the classical one, \eqref{NS_classical}.
\par
The operator \eqref{isom:1-3} being an isomorphism, the space $V_0(\F)$ will be called in the sequel the vorticity space.
This space  is provided with the classical $L^2$--scalar product that is denoted by $(\cdot,\cdot)_{V_0(\F)}$.
\begin{definition}[The projectors $\mathsf P_\F$ and $\mathsf Q_\F$]
\label{def_PQ}
The orthogonal projector onto $V_0(\F)$ in $L^2(\F)$ is denoted by $\mathsf P_{\F}$ and the orthogonal projection 
from $H^1(\F)$ onto $S_0(\F)$ for the semi-norm $(\nabla\cdot,\nabla\cdot)_{\mathbf L^2(\F)}$ is denoted by $\mathsf Q_\F$. 
\end{definition}
%
The projectors $\mathsf P_\F$ and $\mathsf Q_\F$ will play an important role in the proof of Theorem~\ref{main_theo_2}. 
They enjoy the following properties:
%
\begin{prop}
\label{regul_PQ}
For every positive integer $k$, the operators $\mathsf P_\F$ and $\mathsf Q_\F$ map continuously $H^k(\F)$ into $H^k(\F)$.
The mapping $\mathsf P_\F:S_0(\F)\longmapsto V_1(\F)$ is invertible and its inverse is $\mathsf Q_\F:V_1(\F)\longmapsto S_0(\F)$.
\end{prop}
%
The function space $V_1(\F)=V_0(\F)\cap H^1(\F)$ provided with the scalar product:
$$(\omega_1,\omega_2)_{V_1(\F)}=(\nabla\mathsf Q_\F\omega_1,\nabla \mathsf Q_\F\omega_2)_{\mathbf L^2(\F)}\forallt 
\omega_1,\omega_2\in V_1(\F),$$
is densely and continuously included in $V_0(\F)$. 
The dual  of $V_1(\F)$ using $V_0(\F)$ as pivot space is denoted by $V_{-1}(\F)$ so that:
$$V_1(\F)\subset V_0(\F)\subset V_{-1}(\F),$$
is a Gelfand triple of Hilbert spaces. We recall:
%
\begin{theorem}[well posedness of the strong NS equations in vorticity formulation \cite{Lequeurre:2020aa}]
For every vorticity field $\omega^0$ in $V_0(\F)$, 
there exists a unique function $\omega$  in the space:
$$L^2(\mathbb R_+;V_1(\F))\cap C(\mathbb R_+;V_0(\F))\cap H^1(\mathbb R_+;V_{-1}(\F)),$$
satisfying the following Cauchy problem   for every $\theta\in V_1(\F)$:
\begin{subequations}
\begin{alignat}{3}
\label{eq:NS_vorti_main}
\frac{\rm d}{{\rm d}t}(\omega,\theta)_{V_0(\F)}+\nu (\omega,\theta)_{V_1(\F)}
-(\omega\nabla^\perp\psi,\nabla\mathsf Q_{\F}\theta)_{\mathbf L^2(\F)}&=0&\quad&\text{on }\mathbb R_+,\\
\omega(0)&=\omega^0&&\text{in }\F,
\end{alignat}
\end{subequations}
where $\psi$ is the stream function deduced from $\omega$ by means of the Biot-Savart operator.
\end{theorem}
%
\begin{rem}
\label{palenstro}
We emphasize  that in the equation \eqref{eq:NS_vorti_main}, the dissipative term is $
\|\omega\|_{V_1(\F)}^2=\|\nabla\mathsf Q_\F\omega\|_{\mathbf L^2(\F)}^2$ and not $\|\nabla\omega\|_{\mathbf L^2(\F)}^2$ as it could be envisioned.
\end{rem}
An other consequence of Proposition~\ref{prop:iso} is the restatement of the Courant--Fischer min-max principle \eqref{eq:rayleigh} for the eigenvalues 
of the Stokes operator, in term of  stream functions. Thus, for every positive integer $k$, we have:
\begin{equation}
\label{eq:rayleigh_stream}
\lambda_k^\F=\min_{V\in W_k^\F}\max_{{\theta\in V}\atop{\theta\neq 0}}\frac{\|\theta\|_{S_1(\F)}^2}{\|\theta\|_{S_0(\F)}^2},
\end{equation}
where $W_k^\F$ stands the set of all the subspaces of $S_1(\F)$ of dimension $k$. We can defined as well a family 
$\big\{\psi_k^\F,\,k\geqslant 1\big\}$ in $S_1(\F)$ that is an orthonormal Riesz basis in $S_0(\F)$ and orthogonal in $S_1(\F)$ and such that:
$$(\psi_k^\F,\theta)_{S_1(\F)}=\lambda_k^\F (\psi_k^\F,\theta)_{S_0(\F)}\forallt \theta\in S_1(\F).$$
%
 \section{Decomposition of the vorticity space}
\begin{figure}
 \input{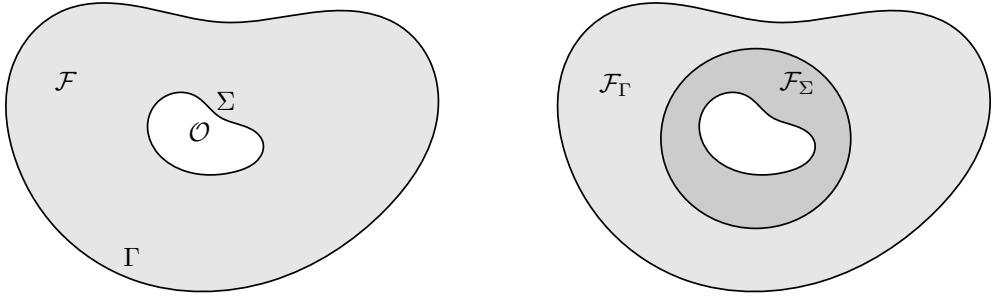}
 \caption{\label{fig1}On the left: The domains and the boundaries; $\G=\F\cup\overline{\B}$. On the right: the partition of $\F$ into $\F_\Sigma$ (a neighborhood of 
 the boundary $\Sigma$) and $\F_\Gamma$.}
 \end{figure}

\begin{figure}
 \input{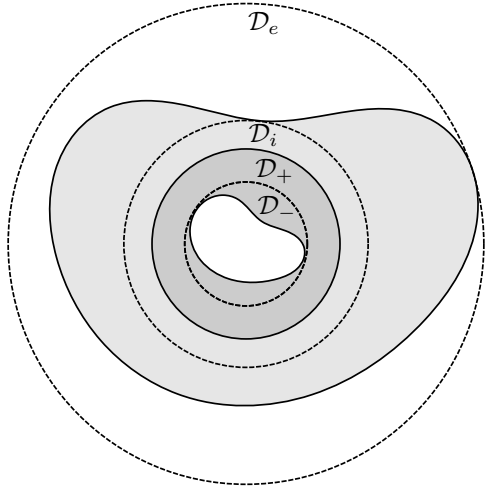}
 \caption{\label{four_concen}The four concentric disks introduced in Hypothesis~\ref{hyp:1}. We recall also the definition of the annuli $\mathcal C_e=\mathcal D_e\setminus\overline{\mathcal D_+}$ and  $\mathcal C_+=\mathcal D_+\setminus\overline{\mathcal D_-}$.}
 \end{figure}

%
\label{SEC:decomp_vorti}
 In this section we consider a fixed perforated  domain $\F=\G\setminus\overline{\B}$ with $\overline{\B}\subset \G$, the sets $\G$ and $\B$ being 
 open, bounded and simply connected, with smooth boundaries denoted respectively by $\Gamma$ and $\Sigma$ (see the left hand side 
 picture of Fig.~\ref{fig1}).  
 We are not yet interested in letting the obstacle $\B$ shrink into a point but we shall assume that, loosely speaking,  
 $\B$ is small enough or far enough from the boundary $\Gamma$ for 
allowing an annulus encircling $\B$ to be included in $\F$. The precise statement of the geometric hypothesizes on the domains is 
rather technical:
 \begin{hypothesis}
 \label{hyp:1}
There exist four concentric disks $\mathcal D_e$, $\mathcal D_i$, $\mathcal D_+$ and $\mathcal D_-$ (see Fig.~\ref{four_concen}) such that:
 $$\overline{\G}\subset \mathcal D_e,\quad \overline{\mathcal D_i}\subset \G,\quad \overline{\mathcal D_+}\subset \mathcal D_i,
 \quad \overline{\mathcal D_-}\subset \mathcal D_+\quad\text{and}\quad \overline{\B}\subset \mathcal D_-.$$
 The radii of the disks $\mathcal D_e, \mathcal D_i, \mathcal D_+$ and $\mathcal D_-$ are denoted respectively by $R_e, R_i, R_+$ and $R_-$ 
 ($R_e$ and $R_i$ are fixed while $R_+$ and $R_-$ are meant to tend to $0$). 
  They are such that:
  \begin{equation}
  \label{rapp_R}
  R_-=R_ee^{-\delta^2},\quad R_+=R_e e^{-\delta}\quad\text{ for some }\delta>\delta_0\text{ with }\delta_0=2+\ln(R_e/R_i).
\end{equation}
 \end{hypothesis}
 %
  In the case where the obstacle $\mathcal O_\varepsilon$ shrinks into a point, the hypothesis above is obviously
 satisfied for $\varepsilon$ small enough providing that:
\begin{equation}
\label{conddde}
\delta\leqslant \sqrt{\ln(R_e/\varepsilon)}.
\end{equation}
 We define the annuli 
  $\mathcal C_e=\mathcal D_e\setminus\overline{\mathcal D_+}$ and  $\mathcal C_+=\mathcal D_+\setminus\overline{\mathcal D_-}$. 
When the obstacle shrinks, the analysis  requires   distinguishing  between the behavior of some quantities near the boundary $\Sigma$ 
and far from this boundary.  To 
this purpose, we introduce also the domains $\F_\Sigma=\mathcal D_+\setminus\overline{\B}$ (a vanishing 
neighborhood of the boundary $\Sigma$ 
in $\F$ where some boundary layer phenomena will take place) and $\F_\Gamma=\F\setminus\overline{\mathcal D_+}$  its supplement in $\F$ (see the right hand side picture of Fig.~\ref{fig1}).
The size of the ``boundary layer'' $\F_\Sigma$ is  given by the identities \eqref{rapp_R}.
%
 %
 %
 %
 \par
Considering the decomposition \eqref{decomp_L2}, we may imagine that the structure of $V_0(\G)$ is tightly related to the structure of $\mathfrak H(\G)$.
The analysis of this latter space is carried out in a series of Lemmas. In the sequel, we shall denote by $f^\star$ a function in $L^2(\G)$ obtained 
by extending by 0 a function $f$ of $L^2(\F)$ or $L^2(\B)$.
 \begin{lemma}
 \label{LEM:dense1}
 Both following assertions hold true:
 \begin{enumerate}
 \item 
  \label{LEM:dense1-1}
 The space $\mathfrak H(\B)$ is the closure in $L^2(\B)$ of the space
 $\{h|_{\B}\,:\,h\in\mathfrak H(\G)\}$.
 \item
  \label{LEM:dense12}
 The space $\mathfrak H_\G(\F)=\{h|_{\F}\,:\,h\in\mathfrak H(\G)\}$
 is closed in $L^2(\G)$. We denote by $\mathfrak H_0(\F)$ its orthogonal complement in $\mathfrak H(\F)$ so that:
\begin{equation}
\label{decomp_ortho_H+}
\mathfrak H(\F)=\mathfrak H_\G(\F)\sumperp\mathfrak H_0(\F).
\end{equation}
 The space $\mathfrak H_0(\F)$ is never reduced to $\{0\}$ and thereby
\begin{equation}
\label{def_h0star}
{\mathfrak H}^\star_0(\F)=\{h^\star\,:\,h\in \mathfrak H_0(\F)\},
\end{equation}
 is a nonempty closed subspace of $V_0(\G)$.
 \end{enumerate}
 \end{lemma}
  \begin{proof}
 Concerning the first assertion, assume that there exists a function $h_\B\in\mathfrak H(\B)$ such that:
  $$\int_{\B}h_\B h_\G\dx=0\qquad\forallt h_\G\in\mathfrak H(\G).$$
  This means that $h^\star_\B$ (i.e. $h_\B$ extended by 0 in $\F$) belongs to $V_0(\G)$ and therefore (using the isometry 
  \eqref{isom:1-3} of Proposition~\ref{isom:1}), there exists a function 
  $\psi\in H^2_0(\G)$ such 
  that $\Delta\psi=h_\B^\star$. But in $\F$ the function $\psi$ is harmonic and on the boundary $\Gamma$ it satisfies $\psi|_{\Gamma}=0$ and 
  $\partial\psi/\partial n|_{\Gamma}=0$. The unique continuation principle for harmonic functions asserts that $\psi=0$ in $\F$.
  Since $\psi$ belongs to $H^2_0(\G)$, we deduce that $\psi|_{\Sigma}=\partial\psi/\partial n|_{\Sigma}=0$ and therefore that $h_\B$ is in 
  $V_0(\B)$ (using again the isometry \eqref{isom:1-3} of Proposition~\ref{isom:1}). By definition $V_0(\B)=\mathfrak H(\B)^\perp$ and hence $h_\B=0$.
  \par
  Addressing the second assertion of the lemma, observe that the convergence of a sequence of functions in $\mathfrak H_\G(\F)$ entails in particular 
  the convergence of the traces of these functions in 
  $H^{-1/2}(\Gamma)$ and therefore of the functions in $\mathfrak H(\G)$. This proves that $\mathfrak H_\G(\F)$ is closed.
  \par
  Let $h_\G$ be a nonconstant function in $\mathfrak H(\G)$ and define $g_0=h_\G|_{\Sigma}$ and $g_1=\partial h_\G/\partial n|_\Sigma$ 
  (notice that $g_0$ cannot be constant because this would entails that $h_\G$ is constant in $\B$ and then also in $\G$). Let $\psi$ be the biharmonic function in $\F$ such that $\psi|_{\Gamma}=0$, $\partial\psi/\partial n|_{\Gamma}=0$ and  $\partial\psi/\partial n|_{\Sigma}=g_1$, 
  $\psi|_{\Sigma}=g_0+c$ with $c$ a constant such that $\int_\Sigma\partial(\Delta\psi)/\partial n\ds=0$ (this condition makes sense since 
  $\partial(\Delta\psi)/\partial n\in H^{-1/2}(\Sigma)$).
 Then, for every $h'_\G\in \mathfrak H(\G)$:
 $$\int_\F \Delta\psi\,h'_\G\,\dx=\int_\Sigma \Big(g_1\,h'_\G-g_0\,\frac{\partial h'_\G}{\partial n}\Big)\ds
 =\int_\B\Big(\Delta h'_\G\,h_\G-h'_\G \Delta h_\G\Big)\dx=0,$$
 and therefore $\Delta\psi$ belongs to $\mathfrak H_0(\F)$. Moreover the function $\psi$ cannot be harmonic in $\F$ because 
its boundary conditions on $\Gamma$ would imply that $\psi=0$ in $\F$ (according to the unique continuation principle) and then that 
$g_0$ is constant. The proof of the lemma is now complete.
  \end{proof}
   We shall establish now some $L^2$--mass concentration properties for harmonic functions. We remind that $R_i$ and $R_e$ are fixed (and depend only on $\G$) 
   while $R_+$ and $R_-$ are meant to tend to 0.
 %
 \begin{lemma}
 \label{LEM:rapp}
Under Hypothesis~\ref{hyp:1}, the following estimate holds:
 $$\sup_{h\in\mathfrak H(\G)}\frac{\|h\|_{L^2(\B)}}{\|h\|_{L^2(\F)}}\leqslant \sqrt{\frac{R_-^2}{R_i^2-R_-^2}}.$$
  \end{lemma}
  \begin{proof}
  In $\mathcal D_i$, the function $h$ can be expanded in polar coordinates as:
  $$h(r,\theta)=a_0+\sum_{n\geqslant 1}r^n\big[a_n\cos(n\theta)+b_n\sin(n\theta)\big],$$
  where $(a_n)_{n\geqslant 0}$ and $(b_n)_{n\geqslant 1}$ are two sequences of real numbers. 
%
Straightforward  computations lead to:
  $$\|h\|_{L^2(\mathcal D_-)}^2=\left(\frac{R_-}{R_i}\right)^2\left[\pi a_0^2 R^2_i+\pi\sum_{n\geqslant 1}\left(\frac{a_n^2+b_n^2}{2n+2}\right) R^{2n+2}_i
  \left(\frac{R_-}{R_i}\right)^{2n}\right]\leqslant \left(\frac{R_-}{R_i}\right)^2\|h\|_{L^2(\mathcal D_i)}^2.$$
  We denote by $\mathcal C_i$ the annulus $\mathcal D_i\setminus\overline{\mathcal D_-}$ and, because of the inclusions $\overline{\B}\subset \mathcal D_-$ 
  and $\mathcal C_i\subset \F$, we have for every $h\in\mathfrak H(\G)$:
  $$\frac{\|h\|^2_{L^2(\B)}}{\|h\|^2_{L^2(\F)}}\leqslant\frac{\|h\|^2_{L^2(\mathcal D_-)}}{\|h\|^2_{L^2(\mathcal C_i)}}=\frac{\|h\|^2_{L^2(\mathcal D_-)}}
  {1-\|h\|^2_{L^2(\mathcal D_-)}}.$$
  The conclusion follows.
  \end{proof}
We introduce the exterior domain $\K=\mathbb R^2\setminus\overline\B$ and the space:
$$\mathfrak H_{\K}(\F)=\Big\{h|_{\F}\,:\,h\in L^2_{\ell oc}(\K),\,h\text{ harmonic in }\K\text{ and }\lim_{|x|\to+\infty} h(x)=0\Big\}.$$
  The asymptotic behavior of the functions in $\mathfrak H_{\K}(\F)$ entails in particular that the flux condition \eqref{nul_flux} is satisfied 
  for every $h\in \mathfrak H_{\K}(\F)$ and every
  smooth Jordan curve $\mathscr C$ included in $\K$.
%
\par
The space $\mathfrak H(\F)\cap H^1(\F)$ is denoted by $\mathfrak H^1(\F)$. Similarly, we define 
\begin{align*}
\mathfrak H_\G^1(\F)=\big\{h|_\F\,:\, h\in\mathfrak H^1(\G)\big\}=\mathfrak H_\G(\F)\cap H^1(\F),\\
\mathfrak H_\K^1(\F)=\big\{h|_\F\,:\, h\in\mathfrak H^1(\K)\big\}=\mathfrak H_\K(\F)\cap H^1(\F),
\end{align*}
where $\mathfrak H^1(\K)=\mathfrak H(\K)\cap H^1_{\ell oc}(\K)$. Notice in particular that every function in $\mathfrak H^1(\K)$ has finite Dirichlet energy i.e. 
$\|\nabla h\|_{\mathbf L^2(\K)}<+\infty$. It is worth observing that the spaces $\mathfrak H^1(\G)$ and $\mathfrak H^1(\K)$ can also be defined by means of 
single layer
potentials. We refer for instance to the book \cite{McLean:2000aa} for details on the single layer potential. Basics on this topic are also 
available in a section of \cite{Munnier:2017aa}. Thus, we denote by 
$\mathsf S_\Gamma:H^{-1/2}(\Gamma)\longrightarrow \mathfrak H^1(\G)$ the single layer potential on $\Gamma$ defined for every $q\in L^2(\Gamma)$ (and extended by density in $H^{-1/2}(\Gamma)$) by:
$$\mathsf S_\Gamma q(x)=\frac{1}{2\pi}\int_\Gamma\ln|x-s|q(s)\ds\forallt x\in \G.$$
The operator $\mathsf S_\Gamma$ is an isomorphism. Define now $\widehat H^{-1/2}(\Sigma)=\big\{q\in H^{-1/2}(\Sigma)\,:\,\int_\Sigma q\ds=0\big\}$ 
and the operator $\mathsf S_\Sigma:\widehat H^{-1/2}(\Sigma)\longrightarrow \mathfrak H^1(\K)$, for every $q\in L^2(\Sigma)\cap \widehat H^{-1/2}(\Sigma)$, by:
$$\mathsf S_\Sigma q(x)=\frac{1}{2\pi}\int_\Sigma\ln|x-s|q(s)\ds\forallt x\in \K.$$
Then, the operator $\mathsf S_\Sigma$ is also an isomorphism. 
%
  \par
  We gave earlier a first decomposition (see \eqref{decomp_ortho_H+}) of the space $\mathfrak H(\F)$. We shall provide below an other one in terms 
  of the spaces $\mathfrak H_\G(\F)$ and $\mathfrak H_\K(\F)$.
 \begin{lemma}
 \label{lem:n1}
%
The space $\mathfrak H(\F)$ admits the following (non-orthogonal) decomposition:
\begin{equation}
\label{eq:decomp_H+}
\mathfrak H(\F)=\mathfrak H_\G(\F)\oplus \mathfrak H_\K(\F).
\end{equation}
Moreover, under   Hypothesis~\ref{hyp:1}, the following equivalent estimates hold for the functions in $\mathfrak H_\K(\F)$:
\begin{subequations}
\begin{equation}
\label{estimhEc}
\inf_{h\in\mathfrak H_\K(\F)}\frac{\|h\|_{L^2(\F_\Sigma)}}{\|h\|_{L^2(\F)}}\geqslant \sqrt{1-\frac{1}{\delta}}
\quad\text{and}
\quad\sup_{h\in\mathfrak H_\K(\F)}\frac{\|h\|_{L^2(\F_\Gamma)}}{\|h\|_{L^2(\F)}}\leqslant \sqrt{\frac{1}{\delta}}.
\end{equation}
Concerning the functions in $\mathfrak H_\G(\F)$, they satisfy the estimates:
\begin{equation}
\label{estimhEc:1}
\sup_{h\in\mathfrak H_\G(\F)}\frac{\|h\|_{L^2(\F_\Sigma)}}{\|h\|_{L^2(\F)}}\leqslant \sqrt{\frac{R_+^2}{R_i^2-R_-^2}}\quad\text{and}
\quad
\inf_{h\in\mathfrak H_\G(\F)}\frac{\|h\|_{L^2(\F_\Gamma)}}{\|h\|_{L^2(\F)}}\geqslant \sqrt{\frac{R_i^2-(R_-^2+R_+^2)}{R_i^2-R_-^2}}
.\end{equation}

\end{subequations}

 \end{lemma}
According to \eqref{eq:decomp_H+}, every function $h\in\mathfrak H(\F)$ can be decomposed into a sum $h_\G+h_\K$ with 
$h_\G\in\mathfrak H_\G(\F)$ and $h_\K\in\mathfrak H_\K(\F)$. Suppose now that $\B$ shrinks into a point (i.e. $\delta$ tends to $+\infty$ or equivalently 
$R_+$ tends to 0). 
Then, \eqref{estimhEc:1} means that the function $h_\G$ concentrates (as far as the $L^2$--norm is concerned) far 
from the boundary $\Sigma$, namely in $\F_\Gamma$ while, according to the estimate \eqref{estimhEc}, 
the function $h_\K$ concentrates in $\mathcal F_\Sigma$, that is to say along the boundary $\Sigma$. 
So, although the decomposition \eqref{eq:decomp_H+} is not orthogonal, it becomes in some sens ``more and more'' orthogonal as $\B$ shrinks.
  \begin{proof} Every function in $\mathfrak H^1(\F)$ is the sum of a single layer potential 
  with a density supported by 
  $\Gamma$ (a function of $\mathfrak H_\G^1(\F)$) and a single layer potential with a density supported by $\Sigma$ (a function 
  of $\mathfrak H_\K^1(\F)$).  
  The identity \eqref{eq:decomp_H+} follows by density.
  \par
  Let now $h$ be in $\mathfrak H_\K(\F)$, $h\neq 0$. 
We begin with the obvious estimate:
\begin{equation}
\label{eq:estimmm}
\frac{\|h\|_{L^2(\F_\Sigma)}}{\|h\|_{L^2(\F_\Gamma)}}\geqslant \frac{\|h\|_{L^2(\mathcal C_+)}}
{\|h\|_{L^2(\mathcal C_e)}}.
\end{equation}
In the exterior domain $\mathbb R^2\setminus\overline{\mathcal D_-}$, the harmonic function $h$ can be expanded in polar coordinates as:
$$h(r,\theta)=\sum_{k\geqslant 1}r^{-k}\big(a_k\cos(k\theta)+b_k\sin(k\theta)\big),$$
where $(a_n)_{n\geqslant 1}$ and $(b_n)_{n\geqslant 1}$ are two sequences of real numbers.
It follows that, 
%
%
 in the annulus $\mathcal C_+$:
\begin{subequations}
\label{estim_crown}
\begin{equation}
\|h\|_{L^2(\mathcal C_+)}^2=\pi (a_1^2+b_1^2)\ln(R_+/R_-)+\sum_{k\geqslant 2}\pi\left(\frac{a_k^2+b_k^2}{2k-2}\right)
\left(\frac{1}{R_+}\right)^{2k-2}\left[\left(\frac{R_+}{R_-}\right)^{2k-2}-1\right],
\end{equation}
while in the annulus $\mathcal C_e$:
\begin{equation}
\|h\|_{L^2(\mathcal C_e)}^2=\pi (a_1^2+b_1^2)\ln(R_e/R_+)+\sum_{k\geqslant 2}\pi\left(\frac{a_k^2+b_k^2}{2k-2}\right)
\left(\frac{1}{R_+}\right)^{2k-2}\left[1-\left(\frac{R_+}{R_e}\right)^{2k-2}\right].
\end{equation}
\end{subequations}
We deduce from both identities \eqref{estim_crown} that:
$$\|h\|_{L^2(\mathcal C_e)}^2\leqslant \max\left\{\frac{\ln (R_+/R_e)}{\ln (R_-/R_+)},(R_+/R_e)^2\right\}\|h\|_{L^2(\mathcal C_+)}^2,$$
and from \eqref{rapp_R} that:
$$(R_+/R_e)^2\leqslant \frac{\ln (R_+/R_e)}{\ln (R_-/R_+)}\leqslant \frac{1}{\delta-1}.$$
%
%
All together, with \eqref{eq:estimmm} and
%
since $\|h\|_{L^2(\F)}^2=\|h\|_{L^2(\F_\Sigma)}^2+\|h\|_{L^2(\F_\Gamma)}^2$, we prove \eqref{estimhEc}. 
\par
Finally, for every $h$   in $\mathfrak H(\G)$, proceeding as in the proof of Lemma~\ref{LEM:rapp}, we easily show that:
$$\frac{\|h\|_{L^2(\F_\Sigma)}}{\|h\|_{L^2(\F)}}\leqslant \frac{\|h\|_{L^2(\mathcal D_+)}}{\|h\|_{L^2(\mathcal C_i)}}
\leqslant \sqrt{\frac{R_+^2}{R_i^2-R_-^2}},$$
and we conclude  the proof.
    \end{proof}
%
Slight modifications
in the proof of the   lemma~\ref{lem:n1} lead to the statement of its counterpart in terms of the $H^1$--norm.
%
 \begin{lemma}
 \label{lem:n2}
The following estimate holds true:
\begin{equation}
\label{inf_norm_H1}
\sup_{h\in\mathfrak H^1(\K)}\frac{\|\nabla h\|_{\mathbf L^2(\mathbb R^2\setminus\overline{\G})}}{\|\nabla h\|_{\mathbf L^2(\F)}}\leqslant
\sqrt{\frac{R_-^2}{R_i^2-R_-^2}}.
\end{equation}
The space $\mathfrak H^1(\F)$ admits the following non-orthogonal decomposition:
\begin{equation}
\label{eq:decomp_H1}
\mathfrak H^1(\F)=\mathfrak H_\G^1(\F)\oplus \mathfrak H_\K^1(\F).
\end{equation}
Moreover, under   Hypothesis~\ref{hyp:1}, the following (equivalent) estimates hold for the functions in $\mathfrak H_\K^1(\F)$:
\begin{subequations}
\label{estimhEcH}
\begin{equation}
\label{estimhEcH1}
\inf_{h\in\mathfrak H^1_\K(\F)}\frac{\|\nabla h\|_{\mathbf L^2(\F_\Sigma)}}{\|\nabla h\|_{\mathbf L^2(\F)}}\geqslant \sqrt{1-\left(\frac{R_-}{R_+}\right)^2}\quad\sup_{h\in\mathfrak H^1_\K(\F)}\frac{\|\nabla h\|_{\mathbf L^2(\F_\Gamma)}}{\|\nabla h\|_{\mathbf L^2(\F)}}\leqslant \left(\frac{R_-}{R_+}\right).
\end{equation}
Regarding the functions in $\mathfrak H_\G^1(\F)$, we have:
\begin{equation}
\label{estimhEcH2}
\sup_{h\in\mathfrak H^1_\G(\F)}\frac{\|h\|_{L^2(\F_\Sigma)}}{\|h\|_{L^2(\F)}}\leqslant \sqrt{\frac{R_+^2}{R_i^2-R_+^2}}.
\end{equation}
\end{subequations}
 \end{lemma}
%
According to the space decompositions \eqref{decomp_ortho_H+} and \eqref{eq:decomp_H+}, every function $h_0$ of $\mathfrak H_0(\F)$ can be decomposed as:
\begin{equation}
\label{p:sum}
h_0=h_\K-\Pi_\G h_\K,
\end{equation}
where $\Pi_\G$ is the orthogonal projection on  $\mathfrak H_\G(\F)$ in $\mathfrak H(\F)$. Notice that in \eqref{p:sum}:
$$\|h_\K\|_{L^2(\F)}^2=\|h_0\|_{L^2(\F)}^2+\|{\Pi_\G} h_\K\|_{L^2(\F)}^2.$$
As $\B$ shrinks, $h_0$ tends to $h_\K$ in \eqref{p:sum}, or equivalently ${\Pi_\G} h_\K$ tends to 0. More precisely, we claim:
\begin{lemma}
\label{lem:estimPGhe}
For every $h_\K\in \mathfrak H_\K(\F)$, under   Hypothesis~\ref{hyp:1}, the following estimate holds:
\begin{equation}
\label{estimPGhe}
\|{\Pi_\G} h_\K\|_{L^2(\F)}\leqslant  2\sqrt{\frac{\delta_0}{\delta}}\| h_\K\|_{L^2(\F)}.
\end{equation}
%
\end{lemma}
\begin{proof}
The $L^2$--norm of $\Pi_\G h_\K$ can be expressed as:
$$\|{\Pi_\G} h_\K\|_{L^2(\F)}=\sup_{\theta\in \mathfrak H_\G(\F)\atop \theta\neq 0}\frac{1}{\|\theta\|_{L^2(\F)}}(h_\K, \theta)_{L^2(\F)},
$$
and we have:
$$\frac{1}{\|\theta\|_{L^2(\F)}}|(h_\K, \theta)_{L^2(\F)}|\leqslant \|h_\K\|_{L^2(\mathcal F_\Sigma)}\frac{\|\theta\|_{L^2(\mathcal F_\Sigma)}}
{\|\theta\|_{L^2(\F)}}
+ \frac{\|h_\K\|_{L^2(\mathcal F_\Gamma)}}{\|h_\K\|_{L^2(\F)}}\frac{\|\theta\|_{L^2(\mathcal F_\Gamma)}}{\|\theta\|_{L^2(\F)}}
\|h_\K\|_{L^2(\F)}.$$
Using the estimate \eqref{estimhEc:1} for the first term in the right hand side and the estimate \eqref{estimhEc} for the second term, we obtain 
that:
$$\|{\Pi_\G} h_\K\|_{L^2(\F)}\leqslant \left(\sqrt{\frac{R_+^2}{R_i^2-R_-^2}}+\sqrt{\frac{1}{\delta}}\right)
\|h_\K\|_{L^2(\F)},$$
and \eqref{estimPGhe} follows, taking into account \eqref{rapp_R}.
\end{proof}
%

%
We can now address the structure of the space $V_0(\G)$.
%
 \begin{lemma}
 \label{LEM:defT}
For every $h_\B$ in $\mathfrak H(\B)$, there exists a unique function
 $\mathsf Th_\B$ in $\mathfrak H_\G(\F)$ such that the function $\omega$ defined in $\G$ by:
 $$\omega|_{\B}=h_\B\qquad\text{and}\qquad \omega|_{\F}=\mathsf Th_\B,$$
 belongs to $V_0(\G)$.  Under Hypothesis~\ref{hyp:1}, the mapping $\mathsf T: h_\B\in \mathfrak H(\B)\longmapsto \mathsf T h_\B\in\mathfrak H_\G(\F)$ is bounded and:
\begin{equation}
\label{inequ_T}
\|\mathsf Th_\B\|_{L^2(\F)}\leqslant \sqrt{\frac{R_-^2}{R_i^2-R_-^2}} \|h_\B\|_{L^2(\B)}\forallt h_\B\in\mathfrak H(\B).
\end{equation}
 It follows in particular that
\begin{equation}
\label{def_W0}
W_0(\G)=\Big\{\omega\in L^2(\G)\,:\, \omega|_{\B}=h_\B,\quad \omega|_{\F}=\mathsf Th_\B,\quad h_\B\in\mathfrak H(\B)\Big\},
\end{equation}
  is a closed subspace of $V_0(\G)$.
 \end{lemma}
 Thus, the space $W_0(\G)$ contains piecewise harmonic functions in $\G$ that are orthogonal in $L^2(\G)$ to the harmonic functions in $\G$.
  \begin{proof}
  Let $h_\B$ be given in $\mathfrak H(\B)$ and define 
  $\psi$ in $H^2_0(\G)$ by setting $\Delta\psi=h_\B$ in $\B$ with $\psi=c$ on $\Sigma$ 
  (a constant that will be fixed later) and $\Delta^2\psi=0$ in $\F$ with $\psi|_{\Gamma}=0$, $\partial\psi/\partial n|_{\Gamma}=0$. 
  Choose the constant $c$ such that $\int_\Sigma\partial(\Delta\psi|_\F)/\partial n\ds=0$ (the normal derivative of $\Delta\psi|_\F$ on $\Sigma$ 
  belongs to $H^{-1/2}(\Sigma)$).  Then define $\mathsf Th_\B$ as the orthogonal projection 
  in $\mathfrak H(\F)$ of $\Delta\psi|_{\F}$ on the subspace $\mathfrak H_\G(\F)$. This proves existence. 
  Uniqueness
   is deduced from the first point of Lemma~\ref{LEM:dense1}.
\par
Notice that the function $\mathsf Th_\B$ can equivalently be defined either by:
$$
\mathsf Th_\B=\argmin\Big\{\|\theta\|_{L^2(\F)}\,:\,\theta\in\mathfrak H(\F),\,\int_{\F}\theta h\dx+\int_{\B}h_\B h\dx=0\quad\forall\,h\in\mathfrak H(\G)\Big\},
$$
or by $\mathsf T h_\B=\Delta\psi_\G|_{\F}$ where:
$$\psi_\G=\argmin\big\{\|\Delta\psi\|_{L^2(\G)}\,:\,\psi\in H^2_0(\G),\,\Delta\psi|_{\B}=h_\B\big\}.$$
To prove the estimate \eqref{inequ_T}, let $h_\B$ be given in $\mathfrak H(\B)$. By definition of the space $\mathfrak H_\G(\F)$, the function 
$\mathsf Th_\B$ (defined in $\F$) can be extended in the whole domain $\G$ in such a way that it belongs to $\mathfrak H(\G)$. 
Keeping the same notation for this extended function, it follows that:
$$\|\mathsf Th_\B\|_{L^2(\F)}^2+\int_{\B}h_\B \mathsf Th_\B\dx=0.$$
Invoking now Lemma~\ref{LEM:rapp}, we obtain that:
$$\|\mathsf Th_\B\|_{L^2(\F)}^2\leqslant \|h_\B\|_{L^2(\B)}\|\mathsf Th_\B\|_{L^2(\B)}\leqslant \|h_\B\|_{L^2(\B)}\sqrt{\frac{R_-^2}{R_i^2-R_-^2}}
\|\mathsf Th_\B\|_{L^2(\F)},$$
and the proof is completed.
  \end{proof}
By extending the functions by $0$ in $\B$ or $\F$, we define $V_0^\star(\F)$ and $V_0^\star(\B)$, two closed subspaces of $V_0(\G)$. These spaces, together 
with ${\mathfrak H}^\star_0(\F)$ (defined in \eqref{def_h0star}) and $W_0(\F)$ (defined in \eqref{def_W0}) enter the decomposition of
$V_0(\G)$.

%
\begin{theorem}
\label{theo:main_decomp_conv}
The vorticity space $V_0(\G)$ admits the following orthogonal decomposition:
\begin{equation}
\label{decomp_V0}
V_0(\G)= V_0^\star(\F)\sumperp  V_0^\star(\B)\sumperp {\mathfrak H}_0^\star(\F)\sumperp W_0(\G).
\end{equation}
Under Hypothesis~\ref{hyp:1} and for every $\omega\in V_0(\G)$:
\begin{equation}
\label{theo:conv_1}
\|\omega-\omega_\F\|_{L^2(\G)}\longrightarrow 0\quad\text{ as $\delta$ tends to }+\infty,
\end{equation}
where $\omega_\F$ is  the orthogonal projection of $\omega$ 
on $V_0^\star(\F)$. Moreover, for every $\delta>4\delta_0$ (what means that $\O$ is small enough), we have the following $L^2$ estimate outside 
the boundary layer $\F_\Sigma$:
\begin{equation}
\label{converge_ldb}
\|\omega-\omega_\F\|_{L^2(\F_\Gamma)}\leqslant \sqrt{\frac{2\delta_0}{\delta-4\delta_0}} \|\omega\|_{L^2(\F)}.
\end{equation}
\end{theorem}
 In \eqref{converge_ldb}, the domain $\F_\Gamma$ depends also on $\delta$ since $\F_\Gamma=\G\setminus\overline{\mathcal D_+}$ 
 and the radius $R_+=R_e e^{-\delta}$ of the disk $\mathcal D_+$ tends to $0$ when $\delta$ goes to $+\infty$. Recall that when the obstacle shrinks (i.e. when 
 $\B=\Beps$), we can choose $\delta=\sqrt{\ln(R_e/\varepsilon)}$ for Hypothesis~\ref{hyp:1} to be satisfied.
 \begin{proof}
Let $\omega$ be in $V_0(\G)$ and define   $\omega_\B=\omega|_{\B}$. The function $\omega_\B$ can 
be decomposed as $\omega_\B^0+\omega_\B^{\mathfrak H}$ with $\omega_\B^0\in V_0(\B)$ and $\omega_\B^{\mathfrak H}\in\mathfrak H(\B)$ (because 
by definition $V_0(\B)$ is the orthogonal complement of $\mathfrak H(\B)$ in $L^2(\B)$). Then extend $\omega_\B$
(keeping the same notation) by setting $\omega_\B=\mathsf T\omega_\B^{\mathfrak H}$ in $\F$ and notice that $\omega_\B$, this extended function, belongs
 to $V_0(\G)$. Introduce 
now $\omega_\F=\omega-\omega_\B$ in $\G$. This function is in $V_0(\G)$ (like $\omega$ and $\omega_\B$) and equal to zero in $\B$. Therefore, its restriction to $\F$  belongs to $(\mathfrak H_\G(\F))^\perp=V_0(\F)\oplus \mathfrak H_0(\F)$ and the proof of the identity \eqref{decomp_V0} 
is completed.
\par
Consider again a function $\omega$ in $V_0(\G)$. According to \eqref{decomp_V0}, 
$\omega$ can be decomposed into the orthogonal sum:
\begin{equation}
\label{decom_omega}
\omega=\omega_\F+\omega_\B+\omega_{\mathfrak H}+\omega_W,
\end{equation}
with $\omega_\F\in V_0^\star(\F)$, $\omega_\B\in V_0^\star(\B)$, $\omega_{\mathfrak H}\in\mathfrak H^\star_0(\F)$ and $\omega_W\in 
W_0(\G)$.
In this sum, every term depends on $\delta$. Notice first that $\omega|_{\B}=\omega_\B|_{\B}+\omega_W|_{\B}$ and this sum is orthogonal in $L^2(\B)$. The 
dominated convergence theorem yields the convergence toward $0$ of $\omega_\B|_\B$ and $\omega_W|_{\B}$ 
in $L^2(\B)$. The latter convergence combined with Lemma~\ref{LEM:defT} yields   
the convergence of $\omega_W$ toward $0$ in $L^2(\G)$. Let us turn our attention now to the term $\omega_{\mathfrak H}$ in \eqref{decom_omega}.
By definition of the orthogonal projection, we have:
\begin{equation}
\|\omega_{\mathfrak H}|_\F\|_{L^2(\F)}=\sup_{h_0\in \mathfrak H_0(\F)\atop h_0\neq 0}\frac{1}{\|h_0\|_{L^2(\F)}}(\omega, h_0)_{L^2(\F)}.
\end{equation}
For every $h_0\in\mathfrak H_0(\F)$  decomposed as in \eqref{p:sum}:
\begin{equation}
\label{lalaland}
\frac{1}{\|h_0\|_{L^2(\F)}}\int_\F\omega h_0\dx=\frac{\|h_\K\|_{L^2(\F)}}{\|h_0\|_{L^2(\F)}}\left[\int_{\F}\omega  \frac{h_\K}{\|h_\K\|_{L^2(\F)}}\dx-
\int_{\F}\omega \frac{{\Pi_\G}h_\K}{\|h_\K\|_{L^2(\F)}}\dx\right].
\end{equation}
However, in the decomposition \eqref{p:sum}, by definition of the space $\mathfrak H_\G(\F)$,
the function ${\Pi_\G} h_\K$ can be supposed to be in  $\mathfrak H(\G)$ (we keep the same notation). From this observation, we deduce that in \eqref{lalaland}:
\begin{equation}
-\int_{\F}\omega \frac{{\Pi_\G}h_\K}{\|h_\K\|_{L^2(\F)}}\dx=
\int_{\B}\omega \frac{{\Pi_\G}h_\K}{\|h_\K\|_{L^2(\F)}}\dx.
\end{equation}
According to \eqref{estimPGhe}, we have when $\delta>4\delta_0$:
\begin{subequations}
\label{estim:plein}
\begin{equation}
\label{estim:plein:1}
\frac{\|h_\K\|_{L^2(\F)}}{\|h_0\|_{L^2(\F)}}\leqslant \sqrt{\frac{\delta}{\delta-4\delta_0}}.
\end{equation}
On the other hand, according to \eqref{estimhEc}:
\begin{equation}
\left|\int_\F\omega  \frac{h_\K}{\|h_\K\|_{L^2(\F)}}\dx\right|\leqslant \|\omega\|_{L^2(\F_\Gamma)}\sqrt{\frac{1}{\delta}}
+\|\omega\|_{L^2(\mathcal F_\Sigma)},
\end{equation}
and the second term tends to $0$ according to the dominated convergence theorem. Finally, considering the last term in \eqref{lalaland}:
\begin{equation}
\left|\int_{\B}\omega \frac{{\Pi_\G}h_\K}{\|h_\K\|_{L^2(\F)}}\dx\right|\leqslant \|\omega\|_{L^2(\B)}\frac{\|{\Pi_\G} h_\K\|_{L^2(\B)}}
{\|{\Pi_\G} h_\K\|_{L^2(\F)}}\frac{\|{\Pi_\G} h_\K\|_{L^2(\F)}}{\|h_\K\|_{L^2(\F)}},
\end{equation}
and according to Lemma~\ref{LEM:rapp} and Lemma~\ref{lem:estimPGhe}:
\begin{equation}
\frac{\|{\Pi_\G} h_\K\|_{L^2(\B)}}
{\|{\Pi_\G} h_\K\|_{L^2(\F)}}\leqslant\sqrt{\frac{R_-^2}{R_i^2-R_-^2}}
\qquad\text{and}\qquad \frac{\|{\Pi_\G} h_\K\|_{L^2(\F)}}{\|h_\K\|_{L^2(\F)}}\leqslant  2\sqrt{\frac{\delta_0}{\delta}}.
\end{equation}
\end{subequations}
Using the estimates \eqref{estim:plein} in the equality \eqref{lalaland}, we conclude the proof of the first convergence result \eqref{theo:conv_1}. 
Let us address now the estimate \eqref{converge_ldb}. Considering back the decomposition \eqref{decom_omega}
, we have:
\begin{equation}
\label{eq:omega-omega+}
(\omega-\omega_\F)|_{\F_\Gamma}=\omega_{\mathfrak H}|_{\mathcal F_\Gamma}+\omega_W|_{\mathcal F_\Gamma}.
\end{equation}
According to the decomposition \eqref{p:sum}, there exists $\omega_\K\in h_\K$ such that:
\begin{equation}
\label{eq:decomp_omega+}
\omega_{\mathfrak H}=\omega_\K-{\Pi_\G}\omega_\K,
\end{equation}
whence we deduce that:
$$\|\omega_{\mathfrak H}\|_{L^2(\F_\Gamma)}
\leqslant \|\omega_\K\|_{L^2(\F_\Gamma)}+\|{\Pi_\G}\omega_\K\|_{L^2(\F)}.$$
Using \eqref{estimhEc} for the first term in the right hand side and \eqref{estimPGhe} for the second, we get :
\begin{subequations}
\label{deux_belles}
\begin{equation}
\|\omega_{\mathfrak H}\|_{L^2(\F_\Gamma)}\leqslant 3\sqrt{\frac{\delta_0}{\delta}}\|\omega_\K\|_{L^2(\F)}.
\end{equation}
On the other hand, proceeding as for  \eqref{estim:plein:1}, we obtain that:
\begin{equation}
\|\omega_\K\|_{L^2(\F)}\leqslant \sqrt{\frac{\delta}{\delta-4\delta_0}}\|\omega_{\mathfrak H}\|_{L^2(\F)}
\leqslant \sqrt{\frac{\delta}{\delta-4\delta_0}}\|\omega \|_{L^2(\F)}.
\end{equation}
\end{subequations}
Combining both estimates \eqref{deux_belles} yields:
\begin{subequations}
\label{deux_der}
\begin{equation}
\|\omega_{\mathfrak H}\|_{L^2(\F_\Gamma)}\leqslant3\sqrt{\frac{\delta_0}{\delta-4\delta_0}} \|\omega\|_{L^2(\F)}.
\end{equation}
Going back to \eqref{eq:omega-omega+} and recalling the definition \eqref{def_W0} of $W_0(\G)$, it comes:
$$\|\omega_W\|_{L^2(\F)}=\|\mathsf T(\omega_W|_\B)\|_{L^2(\F)}.$$
%
We can then apply Lemma~\ref{LEM:defT} to obtain:
\begin{equation}
\|\omega_W\|_{L^2(\F)}\leqslant \sqrt{\frac{R_-^2}{R_i^2-R_-^2}} \|\omega_W|_\B\|_{L^2(\B)}
\leqslant \sqrt{\frac{R_-^2}{R_i^2-R_-^2}} \|\omega\|_{L^2(\F)}.
\end{equation}
\end{subequations}
Using both estimates \eqref{deux_der} in \eqref{eq:omega-omega+}, we obtain \eqref{converge_ldb} and complete the proof.
\end{proof}
%
The last lemma of this section is the cornerstone of the proof of Theorem~\ref{main_theo_2}. It concerns also the behavior of 
harmonic functions. Thus, let $\psi$ be in $S_0(\G)$ (this space is defined in \eqref{def_S_spaces}). The function $\psi|_{\F}$ can be decomposed as 
\begin{equation}
\label{decomp_psieps}
\psi|_{\F}=\mathsf Q_{\F}\psi+h_\F,
\end{equation}
where $h_\F$ belongs to $\mathfrak H^1(\F)$ and the projector $\mathsf Q_\F$ is introduced in Definition~\ref{def_PQ}.  We shall now prove that when the domain $\B$ shrinks (or 
more precisely when $R_-$ tends to 0), the $H^1$--norm in $\F_\Gamma$ of the harmonic function $h_\F$ tends to $0$. 
%
\begin{lemma}
\label{lem36}
Under   Hypothesis~\ref{hyp:1}, there exists a constant $\mathbf c_{[\G]}$ such that, for every $\psi\in S_0(\G)$:
$$\|\nabla (\psi-\mathsf Q_{\F}\psi)\|_{\mathbf L^2(\F_\Gamma)}\leqslant  \mathbf c_{[\G]}\sqrt{\frac{R_-^2}{R_i^2-R_-^2}}\| \psi\|_{S_0(\G)}.$$
\end{lemma}
%
\begin{proof}
According to the identity \eqref{eq:decomp_H1} in Lemma~\ref{lem:n2}, 
the harmonic function $h_\F$ in \eqref{decomp_psieps} can be  decomposed   as:
\begin{equation}
\label{decomp_heps}
h_\F=h_\G+h_{\K}.
\end{equation}
Forming the scalar 
product of \eqref{decomp_heps} with $h_\K$ in $H^1(\F)$, we deduce that:
$$\|\nabla h_\K\|_{\mathbf L^2(\F)}^2
\leqslant \|\nabla h_\F\|_{\mathbf L^2(\F)}\|\nabla h_\K\|_{\mathbf L^2(\F)}
+\|\nabla h_{\G}\|_{\mathbf L^2(\F_\Gamma)}\|\nabla h_{\K}\|_{\mathbf L^2(\F_\Gamma)}
+\|\nabla h_{\G}\|_{\mathbf L^2(\F_\Sigma)}\|\nabla h_{\K}\|_{\mathbf L^2(\F_\Sigma)}.$$
We use now the estimates \eqref{estimhEcH}  to obtain:
$$\|\nabla h_\K\|_{\mathbf L^2(\F)}\leqslant \left( \left(\frac{R_-}{R_+}\right)+
 \sqrt{\frac{R_+^2}{R_i^2-R_+^2}}\right)\|\nabla h_{\G}\|_{\mathbf L^2(\F)}+\|\nabla h_\F\|_{\mathbf L^2(\F)}.
$$
The very same estimate holds true inverting the roles played by $h_\G$ and $h_\K$, whence:
\begin{equation}
\label{eq_estim_hG}
\|\nabla h_\K\|_{\mathbf L^2(\F)}\leqslant \left(1-\left(\frac{R_-}{R_+}\right)-
 \sqrt{\frac{R_+^2}{R_i^2-R_+^2}}\right)^{-1}\|\nabla h_\F\|_{\mathbf L^2(\F)}\leqslant \mathbf c\|\nabla h_\F\|_{\mathbf L^2(\F)},
 \end{equation}
the second inequality resulting from \eqref{rapp_R}.
The decomposition \eqref{decomp_psieps} being orthogonal, $\|\nabla h_\F\|_{\mathbf L^2(\F)}
\leqslant \|\psi\|_{S_0(\G)}$ and therefore there exists a constant $\mathbf c$ such that:
\begin{equation}
\label{unif_bound}
\|\nabla h_\G\|_{\mathbf L^2(\F)}\leqslant\mathbf c \|\psi\|_{S_0(\G)}
\quad\text{ and }\quad \|\nabla h_\K\|_{\mathbf L^2(\F)}\leqslant\mathbf c \|\psi\|_{S_0(\G)}.
\end{equation}
The combination of the  estimate \eqref{inf_norm_H1} with \eqref{unif_bound} yields:
\begin{equation}
\label{eq:egal1}
{\|\nabla h_\K\|_{\mathbf L^2(\mathbb R^2\setminus\overline{\G})}}\leqslant 
\sqrt{\frac{R_-^2}{R_i^2-R_-^2}}{\|\nabla h_\K\|_{\mathbf L^2(\F)}}\leqslant \mathbf c
\sqrt{\frac{R_-^2}{R_i^2-R_-^2}}\|\psi\|_{S_0(\G)}.
\end{equation}
Applying $\tr_\Gamma$ (the trace operator on $\Gamma$ valued in $H^{1/2}(\Gamma)$) to the identity \eqref{decomp_psieps}, taking 
into account \eqref{decomp_heps}, we obtain:
\begin{equation}
\label{eq:egal2}
0=\tr_{\Gamma} h_\F=\tr_{\Gamma} h_\G
+\tr_{\Gamma} h_\K.
\end{equation}
Let us recall now  some elementary results of potential theory (we refer again to the book \cite{McLean:2000aa} or to 
the dedicated section in \cite{Munnier:2017aa}).
The flux condition $\int_\Gamma \partial h_\K/\partial n\ds=0$ entails that the trace of $h_\K$ on $\Gamma$ 
belongs to the following subspace of $H^{1/2}(\Gamma)$:
$$\widehat H^{1/2}(\Gamma)=\Big\{\gamma\in H^{1/2}(\Gamma)\,:\,\int_\Gamma\gamma\,\mathfrak e_\Gamma\ds =0\Big\},$$
 where $\mathfrak e_\Gamma$ stands for the equilibrium 
density of $\Gamma$. For any $\gamma\in \widehat H^{1/2}(\Gamma)$, there exists 
a unique function $h_\gamma$ harmonic in $\mathbb R^2\setminus\G$ such that $\tr_\Gamma h_\gamma=\gamma$ 
and $\|\nabla h_\gamma\|_{\mathbf L^2(\mathbb R^2\setminus \G)}<+\infty$. In $\widehat H^{1/2}(\Gamma)$, the norm:
$$\|\gamma\|_{\widehat H^{1/2}(\Gamma)}=\|\nabla h_\gamma\|_{\mathbf L^2(\mathbb R^2\setminus \G)},$$
is equivalent to the usual norm of $H^{1/2}(\Gamma)$.
\par
Combining \eqref{eq:egal1} and \eqref{eq:egal2}, we deduce first that:
$$
\|\tr_{\Gamma} h_\G\|_{H^{1/2}(\Gamma)} \leqslant \mathbf c
\sqrt{\frac{R_-^2}{R_i^2-R_-^2}}\|\psi\|_{S_0(\G)},
$$
and next (considering the function $h_\G$ as defined in the whole domain $\G$) that:
\begin{subequations}
\label{eq:both}
\begin{equation}
\|\nabla h_\G\|_{\mathbf L^2(\G)}\leqslant \mathbf c_{[\G]}\|\tr_{\Gamma} h_\G\|_{H^{1/2}(\Gamma)}
\leqslant \mathbf c_{[\G]}\sqrt{\frac{R_-^2}{R_i^2-R_-^2}}\|\psi\|_{S_0(\G)}.
\end{equation}
The second estimate in \eqref{estimhEcH1} together with \eqref{unif_bound} leads to:
\begin{equation}
\|\nabla h_\K\|_{\mathbf L^2(\F_\Gamma)}\leqslant\mathbf c \left(\frac{R_-}{R_+}\right)\|\psi\|_{S_0(\G)}.
\end{equation}
\end{subequations}
Using both estimates \eqref{eq:both} in the identity \eqref{decomp_heps}, we conclude the proof of the lemma.
\end{proof}
%
\section{Proof of Theorem~\ref{main:theo} and Corollary~\ref{cor:1}}
\label{SEC:proof_of_theo1}
\begin{proof}[Proof of of Theorem~\ref{main:theo}]
For every positive integer $k$, define $W_k^\G$ the set of the subspaces of $S_1(\G)$ of dimension $k$. Similarly $W^\Feps_k$ stands for the set of the subspaces 
of $S_1(\Feps)$ of dimension $k$ (the spaces $S_0$ and $S_1$ are defined in \eqref{def_S_spaces}). Then, 
the Courant--Fischer min-max principle \eqref{eq:rayleigh_stream} for the eigenvalues of the Stokes operator  read 
as follows for the domains $\G$ and $\Feps$ respectively:
$$\lambda_k^\G=\min_{V\in W_k^\G}\max_{{\theta\in V}\atop{\theta\neq 0}}\frac{\|\theta\|_{S_1(\G)}^2}{\|\theta\|_{S_0(\G)}^2}
\qquad\text{and}\qquad
\lambda_k^\Feps=\min_{V\in W_k^\Feps}\max_{{\theta\in V}\atop{\theta\neq 0}}\frac{\|\theta\|_{S_1(\Feps)}^2}{\|\theta\|_{S_0(\Feps)}^2}.$$
Notice that every function of $S_0(\Feps)$ (or $S_1(\Feps)$) can be seen as a function in $S_0(\G)$ (or $S_1(\G)$) with the same 
norm once 
extended by the suitable constant inside $\Beps$. We can then consider that $S_0(\Feps)\subset S_0(\G)$ and $S_1(\Feps)\subset S_1(\G)$.
From the inclusion $W_k^\Feps\subset W_k^\G$ we deduce straightforwardly that $\lambda_k^\G\leqslant \lambda_k^\Feps$. 
\par
Denote by $\big\{\psi_1^\G,\ldots,\psi_k^\G\big\}$ an orthonormal family in $S_0(\G)$ (and orthogonal in $S_1(\G)$) made of the $k$ first eigenfunctions
of the Stokes operator and let $\mathcal W_k^\G$ be the subspace spanned by the stream functions $\psi_j^\G$. We denote by $\Pi_\Feps$ the orthogonal projection 
from $S_1(\G)$ onto $S_1(\Feps)$ and $\Pi_\Feps^\perp=\rm{Id}-\Pi_\Feps$.
From the convergence result \eqref{theo:conv_1} of Theorem~\ref{theo:main_decomp_conv} and \eqref{isom:1-3} of Proposition~\ref{prop:iso}, 
we deduce that:
$$\eta_k(\varepsilon)=\max_{\theta\in \mathcal W_k^\G\atop \theta\neq 0}\frac{\|\Pi_\Feps^\perp\theta\|_{S_1(\G)}}{\|\theta\|_{S_1(\G)}} \longrightarrow0\quad\text{as}\quad\varepsilon\longrightarrow 0.$$
Considering now the norm of $S_0(\G)$ we have also:
$$\max_{\theta\in \mathcal W_k^\G\atop \theta\neq 0}\frac{\|\Pi_\Feps^\perp\theta\|_{S_0(\G)}}{\|\theta\|_{S_0(\G)}}
\leqslant \eta_k(\varepsilon)\sqrt{\frac{\lambda_k^\G}{\lambda_1^\G}}.$$
Then, by direct computation, we show that for $\eta_k(\varepsilon)$ small enough:
$$\max_{\theta\in \mathcal W_k^\G\atop \theta\neq 0}\left|\frac{\|\theta\|_{S_1(\G)}^2}{\|\theta\|_{S_0(\G)}^2}-\frac{\|\Pi_\Feps\theta\|_{S_1(\G)}^2}{\|\Pi_\Feps\theta\|_{S_0(\G)}^2}\right|
\leqslant \mathbf c_{[\lambda_1^\G,\lambda_k^\G]}\eta_k(\varepsilon),$$
where $\mathbf c_{[\lambda_1^\G,\lambda_k^\G]}$ is a positive constant depending on $\lambda_1^\G$ and $\lambda_k^\G$ only. 
For $\eta_k(\varepsilon)$ small enough, the family $\big\{\Pi_\Feps\psi^\G_1,\ldots,\Pi_\Feps\psi^\G_k\big\}$ is free and spanned a subspace of $S_1(\Feps)$  of 
dimension $k$. If follows that:
$$\lambda_k^\Feps\leqslant \max_{{\theta\in \mathcal W_k^\G}\atop{\theta\neq 0}}\frac{\|\Pi_\Feps\theta\|_{S_1(\G)}^2}{\|\Pi_\Feps\theta\|_{S_0(\G)}^2}
\leqslant \max_{{\theta\in \mathcal W_k^\G}\atop{\theta\neq 0}} \frac{\|\theta\|_{S_1(\G)}^2}{\|\theta\|_{S_0(\G)}^2}
+\max_{\theta\in \mathcal W_k^\G\atop \theta\neq 0}\left|\frac{\|\theta\|_{S_1(\G)}^2}{\|\theta\|_{S_0(\G)}^2}-\frac{\|\Pi_\Feps\theta\|_{S_1(\G)}^2}{\|\Pi_\Feps\theta\|_{S_0(\G)}^2}\right|\leqslant \lambda_k^\G+ \mathbf c_{[\lambda_1^\G,\lambda_k^\G]}\eta_k(\varepsilon),$$
and the proof of \eqref{conv_eigen} is completed. 
\par
Let us address  the result \eqref{conv_eigen:2} about the convergence of the eigenspaces. We  consider again a Riesz orthonormal basis 
$\big\{\psi_j^\G,\,j\geqslant 1\big\}$ of $S_0(\G)$ made of eigenfunctions of the Stokes operator in $\G$ (in stream function formulation). Similarly, 
for every $\varepsilon$, we introduce  $\big\{\psi_j^\Feps,\,j\geqslant 1\big\}$ a Riesz orthonormal basis in $S_0(\Feps)$ 
made of eigenfunctions of the Stokes operator in $\Feps$.
Let a positive integer $k$ be given and let $m\notin  I _k^\G$ (recall that $ I _k^\G$ is the set of all the indices $j$ such that 
$\lambda_j^\G=\lambda_k^\G$). Then:
$$\big(\psi_k^\Feps,\psi_m^\G\big)_{S_0(\G)}=\frac{1}{\lambda_m^\G}\big(\psi_k^\Feps,\Pi_\Feps\psi_m^\G\big)_{S_1(\G)}
=\frac{\lambda_k^\Feps}{\lambda_m^\G}\big(\psi_k^\Feps,\Pi_\Feps\psi_m^\G\big)_{S_0(\G)}.$$
It follows that:
\begin{equation}
\label{eq:first_step:1}
\Big(\psi_k^\Feps,\left(1-\frac{\lambda_k^\Feps}{\lambda_m^\G}\right)\psi_m^\G+\frac{\lambda_k^\Feps}{\lambda_m^\G}
\Pi_\Feps^\perp\psi_m^\G\Big)_{S_0(\G)}=0.
\end{equation}
According to \eqref{conv_eigen}, for $\varepsilon$ small enough, the eigenvalue 
$\lambda_k^\Feps$ is closed to $\lambda_k^\G$ and therefore $\lambda_k^\Feps\neq \lambda_m^\G$.
We deduce first that:
\begin{equation}
\label{ident:lambda}
\big(\psi_k^\Feps,\psi_m^\G\big)_{S_0(\G)}=\frac{\lambda_k^\Feps}{\lambda_k^\Feps-\lambda_m^\G}\big(\psi_k^\Feps,\Pi_\Feps^\perp\psi_m^\G\big)_{S_0(\G)},
\end{equation}
and then, summing over all the indices $m\notin  I _k^\G$:
$$\sum_{m\notin  I _k^\G}\big(\psi_k^\Feps,\psi_m^\G\big)_{S_0(\G)}^2= \big(\lambda_k^\Feps\big)^2\sum_{m\notin I_k^\G}
\frac{\big(\psi_k^\Feps,\Pi_\Feps^\perp\psi_m^\G\big)_{S_0(\G)}^2}{\big(\lambda_k^\Feps-\lambda_m^\G\big)^2}.$$
It is known (see for instance \cite{Ilyin:2009aa}) that $\lambda_m^\G=\mathscr O(m)$  as $m\to +\infty$. On the other hand, for every $m$:
$$\big|\big(\psi_k^\Feps,\Pi_\Feps^\perp\psi_m^\G\big)_{S_0(\G)}^2\big|\leqslant \frac{1}{\lambda_1^\G}\|\Pi_\Feps^\perp\psi_m^\G\|_{S_1(\G)}^2,$$
and this quantity tends to $0$ along with $\varepsilon$ according to the convergence result 
\eqref{theo:conv_1} of Theorem~\ref{theo:main_decomp_conv}. The dominated convergence Theorem ensures next that:
\begin{equation}
\label{tend0:1}
\sum_{m\notin  I _k^\G}\big(\psi_k^\Feps,\psi_m^\G\big)_{S_0(\G)}^2\longrightarrow 0\quad\text{ as }\quad\varepsilon\longrightarrow 0.
\end{equation}
Define now $\underline{\varLambda}_k^\Feps=\bigoplus_{j\in I_k^\G}\varLambda_j^{\Feps}$ (recall that $\varLambda_j^{\Feps}$ 
is the eigenspace associated to the eigenvalue $\lambda_j^\Feps$). Then, for every $\theta\in S_0(\G)$:
$$\Pi_{\varLambda_k^\G}^\perp\Pi_{\underline{\varLambda}_k^\Feps}\theta
=\sum_{m\notin  I _k^\G}\left(\sum_{j\in  I _k^\G}(\theta,\psi_j^\Feps)_{S_0(\G)}(\psi_j^\Feps,\psi_m^\G)_{S_0(\G)}\right)\psi_m^\G,$$
whence 
we deduce that:
\begin{equation}
\label{eq:ZEDFR}
\|\Pi_{\varLambda_k^\G}^\perp\Pi_{\underline{\varLambda}_k^\Feps}\theta\|^2_{S_0(\G)}
\leqslant \|\theta\|_{S_0(\G)}^2\sum_{j\in  I _k^\G}\left(\sum_{m\notin  I _k^\G}(\psi_j^\Feps,\psi_m^\G)_{S_0(\G)}^2\right),
\end{equation}
and the double sum in the right hand side tends to 0 as $\varepsilon$ goes to 0 according to \eqref{tend0:1}.
\par
Let consider back the identity \eqref{eq:first_step:1}, switching the indices $k$ and $m$:
\begin{equation}
\label{eq:first_step:2}
\Big(\psi_m^\Feps,\left(1-\frac{\lambda_m^\Feps}{\lambda_k^\G}\right)\psi_k^\G+\frac{\lambda_m^\Feps}{\lambda_k^\G}
\Pi_\Feps^\perp\psi_k^\G\Big)_{S_0(\G)}=0.
\end{equation}
Denote by $k^-$ the lowest index in $I_k^\G$ and by $k^+$ the largest index.  Recall that the indice $m$ is assumed not belonging to 
$ I _k^\G$. It means that either $m
\leqslant k^--1$ and we can assume that for $\varepsilon$ small enough 
$\lambda_m^\Feps$ is closed to $\lambda_m^\G$ or $m\geqslant k^++1$ and for every $\varepsilon$, $\lambda_m^\Feps\geqslant \lambda_{k^++1}^\G>\lambda_k^\G$. In either case, for every $\varepsilon$ small enough, $\lambda_m^\Feps\neq \lambda_k^\G$.  
We deduce  that:
$$\big(\psi_m^\Feps,\psi_k^\G\big)_{S_0(\G)}=\frac{\lambda_m^\Feps}{\lambda_m^\Feps-\lambda_k^\G}\big(\psi_m^\Feps,\Pi_\Feps^\perp\psi_k^\G\big)_{S_0(\G)},$$
and then, summing over all the indices $m\notin  I _k^\G$:
$$\sum_{m\notin  I _k^\G}(\psi_k^\G,\psi_m^\Feps)_{S_0(\G)}^2=\sum_{m\notin  I _k^\G}
\left(\frac{\lambda_m^\Feps}{\lambda_m^\Feps-\lambda_k^\G}\right)^2(\Pi_\Feps^\perp\psi_k^\G,\psi_m^\Feps)_{S_0(\G)}^2.$$
In the right hand side, the first term in the sum is uniformly bounded (with respect to $\varepsilon$ and $m$) and for the second,  Parseval's identity yields:
$$\sum_{m\geqslant 1}(\Pi_\Feps^\perp\psi_k^\G,\psi_m^\Feps)_{S_0(\G)}^2=\|\Pi_\Feps^\perp\psi_k^\G\|^2_{S_0(\G)}\leqslant 
\frac{1}{\lambda_1^\G}\|\Pi_\Feps^\perp\psi_k^\G\|^2_{S_1(\G)}.$$
Altogether, we have proved that
$$\sum_{m\notin  I _k^\G}(\psi_k^\G,\psi_m^\Feps)_{S_0(\G)}^2\longrightarrow 0\quad\text{ as }\quad\varepsilon\longrightarrow 0.$$
Noticing now that:
$$\|\Pi_{\varLambda_k^\G} \Pi_{\underline{\varLambda}_k^\Feps}^\perp\theta\|^2_{S_0(\G)}
\leqslant \|\theta\|_{S_0(\G)}^2\sum_{k\in  I _k}\left(\sum_{m\notin  I _k}(\psi_k^\G,\psi_m^\Feps)_{S_0(\G)}^2\right),$$
we deduce with \eqref{eq:ZEDFR} that, for every $\theta\in S_0(\G)$, $\theta\neq 0$:
$$\frac{\|\Pi_{\underline{\varLambda}_k^\Feps}\theta-\Pi_{\varLambda_k^\G}\theta\|_{S_0(\G)}^2}{\|\theta\|_{S_0(\G)}^2}=
\frac{\|\Pi_{\varLambda_k^\G}^\perp\Pi_{\underline{\varLambda}_k^\Feps}\theta\|_{S_0(\G)}^2+\|\Pi_{\varLambda_k^\G} \Pi_{\underline{\varLambda}_k^\Feps}^\perp\theta\|^2_{S_0(\G)}}{\|\theta\|_{S_0(\G)}^2}\longrightarrow 0\quad\text{ as }\varepsilon\longrightarrow 0,$$
and the proof of \eqref{conv_eigen:2}  is completed.
\end{proof}
%
\begin{proof}[Proof of Corollary~\ref{cor:1}]
The semigroup of the Stokes operator in the domain $\G$ reads:
\begin{equation}
\label{def_semigroup_stream}
\mathsf T_\G(t)\theta=\sum_{j\geqslant 1}\big(\theta,\psi_j^\G\big)_{S_0(\G)}e^{-\lambda_j^\G t}\psi_j^\G,\forallt t\geqslant 0
\text{ and }\theta\in S_0(\G).
\end{equation}
Notice that although we use the stream function formulation, we keep the same notation as in \eqref{def_semigroup}. In \eqref{def_semigroup_stream}, we reuse the Riesz orthonormal basis 
$\big\{\psi_j^\G,\,j\geqslant 1\big\}$ of $S_0(\G)$ made of eigenfunctions of the Stokes operator in $\G$ that was
 introduced in the proof of Theorem~\ref{main:theo} above. With similar notation for the Stokes semigroup in the domain $\Feps$, we have, 
 for every $\theta\in S_0(\Feps)$ and every $t\geqslant 0$:
 $$e^{\lambda_1^\G t}\big(\mathsf T_\G(t)\theta-\mathsf T_\Feps(t)\theta)=
 \sum_{j\geqslant 1}e^{(\lambda_1^\G-\lambda_j^\G)t}\Big(\big(\theta,\psi_j^\G\big)_{S_0(\G)}\psi_j^\G-
 \big(\theta,\psi_j^\Feps\big)_{S_0(\G)}e^{(\lambda_j^\G-\lambda_j^\Feps) t}\psi_j^\Feps\Big).$$
 For every positive integer $k$, define the spaces $\overline{\varLambda}_k^\G=\bigoplus_{j=1}^k\varLambda_j^{\G}$ 
 and   $\overline{\varLambda}_k^\Feps=\bigoplus_{j=1}^k\underline{\varLambda}_j^{\Feps}$ where we recall that 
 $\underline{\varLambda}_k^\Feps=\bigoplus_{j\in I_k^\G}\varLambda_j^{\Feps}$. 
 \par
 Let now $\theta$ be fixed in $S_0(\G)$ and any $\zeta>0$ be given. 
Let $N$ be an integer large enough such that:
 $$\|\Pi_{\overline{\varLambda}_N^\G}^\perp \theta\|_{S_0(\G)}\leqslant \zeta.$$
 According to Theorem~\ref{main:theo}, for $\varepsilon$ small enough:
 $$
 \|\Pi_{\overline{\varLambda}_N^\Feps}^\perp \theta-\Pi_{\overline{\varLambda}_N^\G}^\perp \theta\|_{S_0(\G)}
 =\|\Pi_{\overline{\varLambda}_N^\Feps}\theta-\Pi_{\overline{\varLambda}_N^\G} \theta\|_{S_0(\G)}\leqslant\zeta.$$
 Denote by $N^+$ the largest index in $I_N^\G$ and notice now that:
 \begin{multline*}
 \|e^{\lambda_1^\G t}\big(\mathsf T_\G(t)\theta-\mathsf T_\Feps(t)\theta)\|_{S_0(\G)}
 \leqslant \sum_{j=1}^{N^+} \Big\|\big(\theta,\psi_j^\G\big)_{S_0(\G)}\psi_j^\G-
 \big(\theta,\psi_j^\Feps\big)_{S_0(\G)}e^{(\lambda_j^\G-\lambda_j^\Feps) t}\psi_j^\Feps\Big\|_{S_0(\G)}\\
 +\|\Pi_{\overline{\varLambda}_N^\G}^\perp\theta\|_{S_0(\G)}+\|\Pi_{\overline{\varLambda}_N^\Feps}^\perp\theta\|_{S_0(\G)}.
 \end{multline*}
 Invoking again Theorem~\ref{main:theo}, the sum in the right hand side can be made smaller than $\zeta$ assuming that $\varepsilon$ 
 is small enough. It follows that for $N$ large enough and $\varepsilon$ small enough:
 $$\|e^{\lambda_1^\G t}\big(\mathsf T_\G(t)\theta-\mathsf T_\Feps(t)\theta)\|_{S_0(\G)}\leqslant 4\zeta,$$
 which concludes the proof of \eqref{conv_eigen-spaces_1}.
\par
To prove \eqref{conv_eigen-spaces_2}, it suffices to notice that:
$$\sup_{\theta\in K\atop \theta\neq 0}\frac{\|\Pi_{\overline{\varLambda}_N^\G}^\perp \theta\|_{S_0(\G)}}{\|\theta\|_{S_0(\G)}}\longrightarrow 0
\quad\text{as}\quad N\longrightarrow +\infty,$$
and also (according to Theorem~\ref{main:theo}, \eqref{conv_eigen:2}):
$$
 \sup_{\theta\in S_0(\G)\atop \theta\neq 0}\frac{\|\Pi_{\overline{\varLambda}_N^\Feps}^\perp \theta-\Pi_{\overline{\varLambda}_N^\G}^\perp \theta\|_{S_0(\G)}}{\|\theta\|_{S_0(\G)}}\longrightarrow 0
\quad\text{as}\quad N\longrightarrow +\infty,$$
and the conclusion follows. 
\end{proof}
\section{Proof of Theorem~\ref{main_theo_2}}
For every $\varepsilon\in(0,\varepsilon_0)$,  the solution to the Navier--Stokes equations in $\Feps$ 
with initial data $\psi_\varepsilon^0\in S_0(\Feps)$ (as defined in Theorem~\ref{theo_weak_NS})
is denoted by $\psi_\varepsilon$. The function $\psi_\varepsilon$ belongs to the space:
$$H^1(\mathbb R_+;S_{-1}(\Feps))\cap C(\mathbb R_+;S_0(\Feps))\cap L^2(\mathbb R_+;S_1(\Feps)).$$
It satisfies $\psi_\varepsilon(0)=\psi_\varepsilon^0$ and for every $\theta\in S_1(\Feps)$:
\begin{equation}
\label{NS_weak}
\frac{\rm d}{{\rm d}t}(\psi_\varepsilon,\theta)_{S_0(\G)}+\nu(\psi_\varepsilon,\theta)_{S_1(\G)}=
(D^2\theta\nabla^\perp\psi_\varepsilon,\nabla\psi_\varepsilon)_{\mathbf L^2(\G)}\quad\text{ on }\mathbb R_+.
\end{equation}
One easily verifies that for any $\varepsilon\in(0,\varepsilon_0)$ and any $t\geqslant 0$:
\begin{equation}
\label{estim_weak1}
\|\psi_\varepsilon(t)\|_{S_0(\G)}\leqslant \|\psi_\varepsilon^0\|_{S_0(\G)}e^{-\nu\lambda_1^\G t}\quad\text{and}\quad
\int_0^{+\infty}\|\psi_\varepsilon(s)\|_{S_1(\G)}^2\ds\leqslant \frac{1}{2\nu}\|\psi^0_\varepsilon\|_{S_0(\G)}^2.
\end{equation}
Assume now that $\psi_\varepsilon^0\longrightharpoonup\psi^0$ weak in $S_0(\G)$ as $\varepsilon$ goes to $0$.
Then, there exists a function $\psi\in L^\infty(\mathbb R_+;S_0(\G))\cap L^2(\mathbb R_+;S_1(\G))$ such that 
up to a subsequence, 
$$\psi_\varepsilon\longrightharpoonup \psi\text{ weak--}\!\star\text{ in }L^\infty(\mathbb R_+;S_0(\G))
\quad\text{and}\quad
\psi_\varepsilon\longrightharpoonup\psi\text{ weak in }L^2(\mathbb R_+;S_1(\G)).$$
We shall now establish some estimates for the time derivative of $\psi_\varepsilon$ and apply a classical result of compactness 
(see for instance \cite[Theorem 5.1]{Lions:1969aa} in the book of Lions).
From this point on, our proof differs from the one in \cite{Iftimie:2006aa}, 
where the authors deduce compactness from Arzel\`a--Ascoli Theorem. 
\par
Thus, let $\Omega$ be a fixed smooth subdomain included in every $\Feps$ (for $\varepsilon$ small enough) 
such that $\G\setminus\overline\Omega$ 
has a finite number of connected components. As usual, every function $\psi\in S_1(\Omega)$ can be seen 
as a function of $S_1(\G)$ once extended by suitable constants in $\G\setminus\overline\Omega$. For every $\theta\in S_1(\Omega)$, 
we deduce from \eqref{NS_weak} that:
\begin{subequations}
\label{EQ:estim_lady}
\begin{equation}
\left|\frac{\rm d}{{\rm d}t}(\psi_\varepsilon,\theta)_{S_0(\G)}\right|\leqslant \nu\|\psi_\varepsilon\|_{S_1(\G)}\|\theta\|_{S_1(\Omega)}
+\|D^2\theta\|_{L^2(\G)}\|\nabla\psi_\varepsilon\|_{\mathbf L^4(\G)}^2\quad\text{a.e. on }\mathbb R_+.
\end{equation}
On the one hand:
\begin{equation}
\|D^2\theta\|_{L^2(\G)}\leqslant \mathbf c_{[\G]}\|\Delta\theta\|_{L^2(\Omega)}=\mathbf c_{[\G]}\|\theta\|_{S_1(\Omega)}.
\end{equation}
On the other hand, using Ladyzhenskaya's inequality (see \cite[Lemma 3.3, page 291]{Temam:1977aa}), we get:
\begin{equation}
\label{estimL4}
\|\nabla\psi_\varepsilon(t)\|_{\mathbf L^4(\G)}^2\leqslant \mathbf c \|\psi_\varepsilon(t)\|_{S_1(\G)}\|\psi_\varepsilon(t)\|_{S_0(\G)}
\quad\text{for a.e. }t\in\mathbb R_+.
\end{equation}
\end{subequations}
Combining \eqref{estim_weak1} and  \eqref{EQ:estim_lady}, we infer that  $\partial_t(\mathsf Q_\Omega\psi_\varepsilon)$ 
is uniformly bounded in $L^2(\mathbb R_+;S_{-1}(\Omega))$.
By definition of $\mathsf Q_\Omega$ as a projector (see Definition~\ref{def_PQ}), we have:
\begin{subequations}
\label{series_EQ}
\begin{equation}
\|\mathsf Q_\Omega\psi_\varepsilon(t)\|_{S_0(\Omega)}\leqslant \|\nabla\psi_\varepsilon(t)\|_{\mathbf L^2(\Omega)}\leqslant \|\psi_\varepsilon(t)\|_{S_0(\G)}\quad\text{for a.e. }t\in\mathbb R_+.
\end{equation}
On the other hand, according to Proposition~\ref{regul_PQ}, $\mathsf Q_\Omega\psi_\varepsilon(t)$ is in $S_0(\Omega)\cap H^2(\Omega)$ 
for a.e. $t$ and:
\begin{equation}
\|\Delta\mathsf Q_\Omega\psi_\varepsilon(t)\|_{L^2(\Omega)}=\|\Delta\psi_\varepsilon(t)\|_{L^2(\Omega)}
\leqslant \|\psi_\varepsilon(t)\|_{S_1(\G)}\quad\text{for a.e. }t\in\mathbb R_+.
\end{equation}
\end{subequations}
From \eqref{series_EQ} and \eqref{estim_weak1} we deduce that the functions $\mathsf Q_\Omega\psi_\varepsilon$ (for every $\varepsilon\in(0,\varepsilon_0)$) 
are uniformly bounded in 
 $L^2(\mathbb R_+;S_0(\Omega)\cap H^2(\Omega))$ where the Hilbert space $S_0(\Omega)\cap H^2(\Omega)$ is 
provided with the norm: 
$$\big[\|\theta\|_{S_0(\Omega)}^2+\|\Delta\theta\|_{L^2(\Omega)}^2\big]^{1/2}.$$
We can now apply 
the result of compactness  \cite[Theorem 5.1]{Lions:1969aa} with the triple of compactly embedded Hilbert spaces:
$$S_0(\Omega)\cap H^2(\Omega)\subset S_0(\Omega)\subset S_{-1}(\Omega).$$
We deduce that   the family of functions $(\mathsf Q_\Omega\psi_\varepsilon)_{\varepsilon>0}$ is precompact 
in $L^2_{\ell oc}(\mathbb R_+;S_0(\Omega))$ and thereby that, up to a subsequence:
\begin{equation}
\label{strongQpsi}
\mathsf Q_\Omega\psi_\varepsilon\longrightarrow \mathsf Q_\Omega\psi\quad\text{ strong in }L^2_{\ell oc}(\mathbb R_+;S_0(\G)).
\end{equation}
We consider now a geometric configuration as described in the beginning of Section~\ref{SEC:decomp_vorti} with in particular two concentric disks $\mathcal D_+$ and 
$\mathcal D_-$ of radii $R_+$ and $R_-$ given in \eqref{rapp_R} for some $\delta$ large enough such that, for every $\varepsilon$ small enough $\overline{\B}_\varepsilon\subset \mathcal D_-$ (what means that the condition \eqref{conddde} between $\delta$ and 
$\varepsilon$ is satisfied). 
We choose $\Omega=\G\setminus\overline{\mathcal D_-}$ and we denote $\Omega_\Gamma=\Omega\setminus\overline{\mathcal D_+}$.  Let 
us recall also 
the definition of the annulus $\mathcal C_+=\mathcal D_+\setminus\overline{\mathcal D_-}$. 
According to Lemma~\ref{lem36} (in which $\Omega$ plays the role of $\F$, $\Omega_\Gamma$ of $\F_\Gamma$ and $\mathcal C_+$ 
of $\F_\Sigma$), we have:
\begin{equation}
\label{eq_eqtim_psiesp}
\|\nabla(\psi_\varepsilon-\mathsf Q_\Omega\psi_\varepsilon)(t)\|_{\mathbf L^2(\Omega_\Gamma)}
\leqslant \mathbf c_{[\G]}\sqrt{\frac{R_-^2}{R_i^2-R_-^2}}\| \psi_\varepsilon(t)\|_{S_0(\G)}\quad\text{for a.e. }t\in\mathbb R_+.
\end{equation}
The same estimate holds true replacing $\psi_\varepsilon$ by $\psi$. We shall now let $\varepsilon$ tends to $0$ and $\delta$ tends 
to $+\infty$ (or equivalently $R_+$ tends to 0) and prove that $\psi_\varepsilon$ converges strongly to $\psi$ in $L^2_{\ell oc}(\mathbb R_+;S_0(\G))$, i.e. remove the projector $\mathsf Q_\Omega$ in \eqref{strongQpsi}.
On the one hand,  a.e. on $\mathbb R_+$:
\begin{equation}
\label{firdstgsh}
\|\psi-\psi_\varepsilon\|^2_{S_0(\G)}=\|\nabla(\psi_\varepsilon-\psi)\|_{\mathbf L^2(\mathcal D_+)}^2+
\|\nabla(\psi_\varepsilon-\psi)\|_{\mathbf L^2(\Omega_\Gamma)}^2,
\end{equation}
and the first term in the right hand side tends to $0$ in $L^1(\mathbb R_+)$ uniformly with respect to $\varepsilon$ as $\delta$ tends to $+\infty$ according to 
the dominated convergence Theorem and the uniform estimate \eqref{estim_weak1}. Concerning the second term, we write that:
\begin{multline}
\label{plein_de}
\|\nabla(\psi_\varepsilon(t)-\psi(t))\|_{\mathbf L^2(\Omega_\Gamma)}
\leqslant \|\nabla(\psi_\varepsilon(t)-\mathsf Q_\Omega\psi_\varepsilon(t))\|_{\mathbf L^2(\Omega_\Gamma)}+\|\nabla(\mathsf Q_\Omega\psi_\varepsilon(t)-\mathsf Q_\Omega\psi(t))\|_{\mathbf L^2(\Omega_\Gamma)}\\
+\|\nabla(\mathsf Q_\Omega\psi(t)-\psi(t))\|_{\mathbf L^2(\Omega_\Gamma)},
\end{multline}
and therefore, with \eqref{eq_eqtim_psiesp}, for every pair $(\delta,\varepsilon)$ satisfying \eqref{conddde}, we have a.e. 
on $\mathbb R_+$:
$$
\|\nabla(\psi_\varepsilon-\psi)\|_{\mathbf L^2(\Omega_\Gamma)}
\leqslant \mathbf c_{[\G]}\sqrt{\frac{R_-^2}{R_i^2-R_-^2}}\big(\| \psi_\varepsilon\|_{S_0(\G)}+\|\psi\|_{S_0(\G)}\big)\\
+
\|\nabla(\mathsf Q_\Omega\psi_\varepsilon-\mathsf Q_\Omega\psi)\|_{\mathbf L^2(\Omega_\Gamma)}.
$$
We can now choose $\delta$ large enough in such a way that the first term in the right hand side is small in 
$L^2(\mathbb R_+)$ for every $\varepsilon$ satisfying \eqref{conddde}.
Then ($\delta$ being fixed) for every $T>0$ we can make the 
last term in the right hand side also small in $L^2(0,T)$  by choosing $\varepsilon$ small enough, according to \eqref{strongQpsi}. 
%
%
All together with \eqref{firdstgsh}, we can now conclude that (up to a subsequence in $\varepsilon$):
$$\psi_\varepsilon\longrightarrow \psi\quad\text{ strong in }L^2_{\ell oc}(\mathbb R_+;S_0(\G)).$$
It is classical (see for instance \cite[pages 76--77]{Lions:1969aa}) to combine this convergence result with the estimate \eqref{estimL4} and obtain that, for 
every $\theta\in S_1(\G)$:
$$(D^2\theta,\nabla^\perp\psi_\varepsilon,\nabla\psi_\varepsilon)_{\mathbf L^2(\G)}\longrightarrow (D^2\theta,\nabla^\perp\psi,\nabla\psi)_{\mathbf L^2(\G)}\quad
\text{weak  in }L^2_{\ell oc}(\mathbb R_+)\quad\text{ as }\quad \varepsilon\longrightarrow 0.$$
For every $\varepsilon\in(0,\varepsilon_0)$ and a.e. $t\in\mathbb R_+$, denote by $\xi_\varepsilon(t)$ 
the linear form $\partial_t\psi_\varepsilon(t)$ in $S_{-1}(\Feps)$. The same computations as for \eqref{EQ:estim_lady} lead for every $\theta\in S_1(\Feps)$ to:
\begin{equation}
\label{nbhjuio}
\big|\langle\xi_\varepsilon,\theta\rangle_{S_{-1}(\Feps), S_1(\Feps)}\big|\leqslant \nu\|\psi_\varepsilon\|_{S_1(\G)}\|\theta\|_{S_1(\Feps)}
+\|D^2\theta\|_{L^2(\Feps)}\|\nabla\psi_\varepsilon\|_{\mathbf L^4(\G)}^2\quad\text{a.e. on }\mathbb R_+.
\end{equation}
Since $S_1(\Feps)$ is a closed subspace of $S_1(\G)$, we can extend by $0$ every linear form $\xi_\varepsilon(t)$ on $S_1(\G)$ (consider
 an orthogonal supplement of $S_1(\Feps)$ in $S_1(\G)$, which is actually made explicit in Theorem~\ref{theo:main_decomp_conv}). 
 We still denote by $\xi_\varepsilon(t)$ this extended form. From 
\eqref{nbhjuio} and the uniform estimates \eqref{estim_weak1} and \eqref{estimL4}, we deduce that $\xi_\varepsilon$ is uniformly 
bounded with respect to $\varepsilon$ in $L^2(\mathbb R_+;S_{-1}(\G))$ and therefore that there exists $\xi\in L^2(\mathbb R_+;S_{-1}(\G))$ 
such that (up to a subsequence):
\begin{equation}
\xi_\varepsilon\longrightharpoonup\xi\quad\text{weak in }L^2(\mathbb R_+;S_{-1}(\G)).
\end{equation}
Consider again an open set $\Omega=\G\setminus \overline{\mathcal D_-}$  and let  $\theta$ be in $\mathscr D(\mathbb R;S_1(\Omega))$. Then, 
for every $\varepsilon$ small enough $\Omega\subset \Feps$ and:
$$\int_{\mathbb R_+} (\psi_\varepsilon(s),\partial_t\theta(s))_{S_0(\G)}\ds=(\psi_\varepsilon^0,\theta(0))_{S_0(\G)}-\int_{\mathbb R_+} \big\langle\xi_\varepsilon(s), \theta(s)\big\rangle_{S_{-1}(\G),S_1(\G)}\ds.$$
Letting $\varepsilon$ go to $0$, we obtain that:
$$\int_{\mathbb R_+} (\psi(s),\partial_t\theta(s))_{S_0(\G)}\ds=(\psi^0,\theta(0))_{S_0(\G)}-\int_{\mathbb R_+} \big\langle\xi(s), \theta(s)\big\rangle_{S_{-1}(\G),S_1(\G)}\ds.$$
This identity being true for every $\theta\in\mathscr D(\mathbb R_+;S_1(\Omega))$ and for every $\Omega$, we conclude with Theorem~\ref{theo:main_decomp_conv} ($\Omega$ plays the role of $\F$) that it is also true for every $\theta\in\mathscr D(\mathbb R_+;S_1(\G))$
and therefore that $\xi=\partial_t\psi$ and $\psi(0)=\psi^0$. 
\par
Consider back an open set $\Omega=\G\setminus \overline{\mathcal D_-}$  and let  $\theta$ be in $S_1(\Omega)$. Then, we can pass to the limit 
in every term of the equation \eqref{NS_weak} and obtain that $\psi$ belongs to the space:
$$L^2(\mathbb R_+;S_1(\G))\cap C(\mathbb R_+;S_0(\G))
\cap H^1(\mathbb R_+;S_{-1}(\G)),$$
and solves the Cauchy problem:
\begin{alignat*}{3}
\frac{\rm d}{{\rm d}t}(\psi,\theta)_{S_0(\G)}+\nu(\psi,\theta)_{S_1(\G)}&=
(D^2\theta\nabla^\perp\psi,\nabla\psi)_{\mathbf L^2(\G)}&\quad&\text{ on }\mathbb R_+,\\
\psi(0)&=\psi^0&&\text{ in }\G.
\end{alignat*}
This is true for every $\Omega$ and therefore, invoking again Theorem~\ref{theo:main_decomp_conv}, this is true for every $\theta\in S_1(\G)$.
\par
%
We shall prove now the convergence results \eqref{gbhnj}. 
For every $\varepsilon$, the vorticity field $\omega_\varepsilon$ belongs to the space:
$$H^1(\mathbb R_+;V_{-1}(\Feps))\cap C(\mathbb R_+;V_0(\Feps)\cap L^2(\mathbb R_+;V_1(\Feps)),$$
and solves for  every $\theta\in V_1(\Feps)$, the Cauchy problem: 
\begin{subequations}
\begin{alignat}{3}
\label{eq:NS_vorti_main_strong}
\frac{\rm d}{{\rm d}t}(\omega_\varepsilon,\theta)_{V_0(\Feps)}+\nu (\omega_\varepsilon,\theta)_{V_1(\Feps)}
+(\omega_\varepsilon\nabla^\perp\psi_\varepsilon,\nabla\mathsf Q_{\Feps}\theta)_{\mathbf L^2(\Feps)}&=0&\quad&\text{on }\mathbb R_+,\\
\omega_\varepsilon(0)&=\omega^0_\varepsilon&&\text{in }\Feps.
\end{alignat}
\end{subequations}
Let now $\chi$ be in $\mathscr D(\G\setminus\{r\})$. Then, for every $\varepsilon$ small enough, the domain $\B_\varepsilon$ and the support 
$\mathcal U$ of the function $\chi$ are disjoint. We choose $\theta=\theta_\varepsilon=\mathsf P_{\Feps}\chi^2\omega_\varepsilon$ in \eqref{eq:NS_vorti_main_strong} and this function 
is indeed in $V_1(\F_\varepsilon)$ because $\mathsf Q_\Feps\mathsf P_\Feps\chi^2\omega_\varepsilon=\chi^2\omega_\varepsilon$ belongs 
to $S_0(\F_\varepsilon)$. On the one hand, we have:
\begin{subequations}
\label{gtfgt}
\begin{equation}
(\omega_\varepsilon,\theta_\varepsilon)_{V_0(\Feps)}=(\omega_\varepsilon,\chi^2\omega_\varepsilon)_{V_0(\Feps)}=
\|\chi\omega_\varepsilon\|_{L^2(\G)}^2.
\end{equation}
On the other hand:
\begin{align}
\nonumber
(\nabla\mathsf Q_\Feps\mathsf P_\Feps\chi^2\omega_\varepsilon,\nabla\mathsf Q_\Feps\omega_\varepsilon)_{\mathbf L^2(\Feps)}
&=(\nabla\chi^2\omega_\varepsilon,\nabla\omega_\varepsilon)_{\mathbf L^2(\G)}\\
&=\|\nabla(\chi\omega_\varepsilon)\|^2_{L^2(\G)}-\|\omega_\varepsilon\nabla\chi\|^2_{L^2(\G)}.
\end{align}
\end{subequations}
Using \eqref{gtfgt} in \eqref{eq:NS_vorti_main_strong}, we obtain:
\begin{multline}
\label{tyui}
\frac{1}{2}\frac{\rm d}{{\rm d}t}\|\chi\omega_\varepsilon\|_{L^2(\G)}^2+\nu \|\nabla(\chi\omega_\varepsilon)\|^2_{L^2(\G)}=\\
\nu \|\omega_\varepsilon\nabla\chi\|^2_{L^2(\G)}-(\omega_\varepsilon\nabla^\perp\psi_\varepsilon,\chi\omega_\varepsilon\nabla\chi)_{\mathbf L^2(\G)}-
(\chi\omega_\varepsilon\nabla^\perp\psi_\varepsilon, \nabla(\chi\omega_\varepsilon))_{\mathbf L^2(\G)}.
\end{multline}
Applying H\"older's inequality to the second term in the right hand side term leads to:
\begin{subequations}
\label{kolop}
\begin{equation}
|(\omega_\varepsilon\nabla^\perp\psi_\varepsilon,\chi\omega_\varepsilon\nabla\chi)_{\mathbf L^2(\G)}|\leqslant \|\nabla\chi\|_{L^\infty(\G)}
\|\omega_\varepsilon\|_{L^2(\G)}
\|\nabla\psi_\varepsilon\|_{\mathbf L^4(\G)}\|\chi\omega_\varepsilon\|_{L^4(\G)},
\end{equation}
Then, Ladyzhenskaya's inequality yields:
\begin{align}
\|\nabla\psi_\varepsilon\|_{\mathbf L^4(\G)}&\leqslant \mathbf c\|\psi_\varepsilon\|_{S_0(\G)}^{1/2}\|\omega_\varepsilon\|_{L^2(\G)}^{1/2},\\
\|\chi\omega_\varepsilon\|_{L^4(\G)}&\leqslant  \mathbf c\|\chi\omega_\varepsilon\|_{L^2(\G)}^{1/2}
\|\nabla(\chi\omega_\varepsilon)\|_{\mathbf L^2(\G)}^{1/2}.
\end{align}
\end{subequations}
Putting together the estimates \eqref{kolop} and applying Young's inequality, we obtain:
\begin{subequations}
\label{kolop_1}
\begin{multline}
|(\omega_\varepsilon\nabla^\perp\psi_\varepsilon,\chi\omega_\varepsilon\nabla\chi)_{\mathbf L^2(\G)}|\leqslant 
\frac{\nu}{2} \|\nabla(\chi\omega_\varepsilon)\|_{\mathbf L^2(\G)}^2\\
+\mathbf c_{[\nu]}\|\nabla\chi\|_{L^\infty(\G)}^2\|\omega_\varepsilon\|_{L^2(\G)}^2+\mathbf c_{[\nu]}\|\psi_\varepsilon\|_{S_0(\G)}^2
\|\omega_\varepsilon\|_{L^2(\G)}^2
\|\chi\omega_\varepsilon\|_{L^2(\G)}^{2}
\end{multline}
The last term in the right hand side of \eqref{tyui} can be dealt with the same way:
\begin{equation}
|(\chi\omega_\varepsilon\nabla^\perp\psi_\varepsilon, \nabla(\chi\omega_\varepsilon))_{\mathbf L^2(\G)}|\leqslant \frac{\nu}{2} 
\|\nabla(\chi\omega_\varepsilon)\|_{\mathbf L^2(\G)}^2+\mathbf c_{[\nu]}\|\psi_\varepsilon\|_{S_0(\G)}^2\|\omega_\varepsilon\|_{L^2(\G)}^2
\|\chi\omega_\varepsilon\|_{L^2(\G)}^{2}.
\end{equation}
\end{subequations}
Going back to \eqref{tyui} and using both estimates \eqref{kolop_1}, we get:
$$
\frac{\rm d}{{\rm d}t}\|\chi\omega_\varepsilon\|_{L^2(\G)}^2+\Big(\nu\lambda_{\mathcal U}^d-\mathbf c_{[\nu]}\|\psi_\varepsilon\|_{S_0(\G)}^2
\|\omega_\varepsilon\|_{L^2(\G)}^2\Big)
 \|\chi\omega_\varepsilon\|^2_{L^2(\G)}\leqslant \mathbf c_{[\nu]}\|\nabla\chi\|_{L^\infty(\G)}^2\|\omega_\varepsilon\|_{L^2(\G)}^2,
 $$
 where  $\lambda_{\mathcal U}^d$ is the smallest eigenvalue of the Dirichlet Laplacian in the domain $\mathcal U$ (the support 
 of the function $\chi$). We can now apply Grönwall's inequality and thanks to the uniform estimates \eqref{estim_weak1}, we prove first that 
 $\|\chi\omega_\varepsilon\|_{L^2(\G)}$ is uniformly bounded in $L^\infty(\mathbb R_+)$ and next, with \eqref{tyui}, that $\|\nabla(\chi\omega_\varepsilon)\|^2_{L^2(\G)}$ is uniformly bounded in $L^2(\mathbb R_+)$. The convergence results \eqref{gbhnj} follow 
 and the proof is completed.
%

\end{document}